\documentclass[12pt]{article}

\usepackage{amsmath,amsthm,amssymb}

\theoremstyle{plain}
\newtheorem{theorem}{Theorem}
\newtheorem{proposition}{Proposition}[section]
\newtheorem{corollary}{Corollary}[section]
\newtheorem{lemma}{Lemma}[section]
\theoremstyle{remark}
\newtheorem{rem}{\sc Remark}[section]

\makeatletter

\@addtoreset{equation}{section}
\makeatother

\newcommand{\al}{\alpha}
\newcommand{\de}{\delta}
\newcommand{\ep}{\varepsilon}
\newcommand{\lam}{\lambda}
\newcommand{\om}{\omega}

\newcommand{\Ga}{\Gamma}
\newcommand{\La}{\Lambda}

\newcommand{\beq}{\begin{eqnarray*}}
\newcommand{\eeq}{\end{eqnarray*}}
\newcommand{\beqn}{\begin{equation}}
\newcommand{\eeqn}{\end{equation}}

\newcommand{\nin}{\noindent}
\newcommand{\pf}{\noindent {\it Proof. \,}}
\newcommand{\ds}{{\rm   d}}
\newcommand{\as}{{\rm   a}}

\newcommand{\vom}{\varOmega}
\def\v2{\vskip2mm}

\begin{document}

\begin{center}
{\Large   The two-sided exit problem\\
 for a random walk on $\mathbb{Z}$ and having infinite variance II }
\vskip6mm
{K\^ohei UCHIYAMA\footnote{Department of Mathematics,Tokyo Institute of Technology, Japan}} \\
\vskip2mm
\end{center}

\vskip8mm

\begin{abstract} Let $F$ be a distribution function on the integer lattice  $\mathbb{Z}$ and  $S=(S_n)$  the random walk  with step distribution $F$. Suppose $S$ is oscillatory and denote by  $U_{\rm a}(x)$ and $u_{\rm a}(x)$ the renewal function and sequence, respectively,  of  the strictly ascending ladder height process associated with $S$. Putting $A(x) =\int_0^x [1-F(t)-F(-t)] dt$, $H(x)=1-F(x)+F(-x)$ we suppose 
 $$A(x)/\big(xH(x)\big) \to -\infty \quad (x\to\infty).$$ 
Under some additional regularity condition on the positive tail of $F$, we show that
$$u_{\rm a}(x) \sim U_{\rm a}(x)[1-F(x)]/|A(x)|$$
as $x\to\infty$ and uniformly for $0\leq x\leq R\in \mathbb{Z}$, as $R\to\infty$
$$P [  S\;  \mbox{leaves $[0,R]$ on its upper side}\, |\, S_0=x] \, \sim\, c^{-1}A(x)u_{\rm a}(x), $$
where  $c= \sum_{n=1}^\infty P[S_{n}>S_0;\, S_k <  S_0$ for $0<k<n]$ and the  regularity condition is satisfied at least if $S$ is recurrent, $\limsup [1-F(x)]/F(-x)<1$, and  $x[1-F(x)]/L(x) $ ($x\geq 1$) is bounded away from zero and infinity for some slowly varying function  $L$. We also  give  some  asymptotic estimates of  the probability that  $S$ visits $R$ before entering the negative half-line for asymptotically stable walks and obtain asymptotic behaviour of the probability that  $R$ is ever hit by $S$ conditioned to avoid the negative half-line forever.
\end{abstract}



\vskip6mm
\section{Introduction}

 This paper is a continuation of \cite{Uexit}.  We use the same notation as in \cite{Uexit}  which together with the setting of the present work and \cite{Uexit} we present below.  
Let $S= (S_n)_{n=0}^\infty$ be  a  random walk (r.w.)  on the integer lattice $\mathbb{Z}$ with i.i.d.  increments and  an initial point  $S_0$ which is an unspecified integer.   Let $X$ be a generic random variable having the same law as the increment $S_1-S_0$. 
For $x\in \mathbb{Z}$ let $P_x$ denote  the law  of the r.w.\,$S$ started at $x$ and  $E_x$ the expectation by $P_x$;  the subscript $x$ is dropped from $P_x$ and $E_x$ if $x=0$. 
We suppose throughout the paper that $S$ is  irreducible,  {\it oscillating}, and
$\sigma^2:=EX^2=\infty.$
For a subset $B\subset (-\infty,\infty)$ such that $B\cap \mathbb{Z} \neq \emptyset$,  denote by $\sigma_B$ the first time when $S$ visits  $B$ after time zero, namely $\sigma_B= \inf\{n\geq 1: S_n\in B\}$. For simplicity, we write $\sigma_x$ for $\sigma_{\{x\}}$.  
As in \cite{Uexit} we shall be primarily  concerned with the $P_x$ probability of the event
$$\La_R =\{\sigma_{(R,\infty)}< T \},$$
where $\vom = (-\infty,-1]$, $T=\sigma_{\vom}$ and $R$ is a positive integer. 
Denote by  $U_\as(x)$ and $u_\as(x)$ ($V_\ds(x)$ and $v_\ds(x)$)  the renewal function and sequence  of the strictly ascending  (weakly descending) ladder height process associated with $S$. Put 
$$A(x) = \int_0^x[1-F(t) -F(-t)]dt\quad  \mbox{and}\quad H(x)= P_x[|X|>x].$$
 In \cite{Uexit} we observed that  $P_x(\La_R)$  always admits the upper bound  
$$P_x(\La_R) \leq {V_\ds(x)}/{V_\ds(R)} \qquad (x\geq 0)$$
 (cf (2.3) of \cite{Uexit}) and  obtained some sufficient conditions in  order that 
\beqn\label{fun_res}
P_x(\La_R) \sim{V_\ds(x)}/{V_\ds(R)} \quad \mbox{uniformly for $0\leq x<R$ as} \; R\to\infty.
\eeqn
One of them is fulfilled by 
\v2
 {\bf (PRS)} \qquad $A(x)/xH(x) \to \infty$  \quad as \; $x\to \infty$
\v2
\noindent
 so  that (\ref{fun_res}) holds under (PRS),
 while  in \cite{Uexit} we also showed that if
 \v2
 {\bf (NRS)} \qquad $A(x)/xH(x) \to -\infty  \quad \mbox{ as}\quad  x\to \infty,$
\v2\noindent
then $P_x(\La_R)= o\big(V_\ds(x)\big/V_\ds(R)\big)$ uniformly for $0\leq x<\de R$.
 In this paper, we obtain the precise  
asymptotic form of $P_x(\La_R)$ in case (NRS) under some additional regularity condition on the positive tail of $F$ that is satisfied at least when $F$ is in the domain of attraction of a stable law with exponent 1, $EX=0$ and $P[S_n>0] \to 1$. 
This result is accompanied by 
the exact asymptotic form  of $u_\as$. 

The condition (PRS) holds  if and only if the r.w. $S$ is  {\it positively relatively stable} (abbreviated as {\it p.r.s.}) (i.e., there exists a positive sequence $B_n$ such that $S_n/B_n \to1$ in probability) (cf. \cite[p.1478]{KM/ifp}, \cite{M}, \cite{R3}). We  shall say $F$ to be  p.r.s. or n.r.s. according as (PRS) or (NRS) holds. Similarly we shall say  $F$ to be recurrent (or transient) if so is the r.w. $S$.
 If $F$ is p.r.s. (n.r.s.)  both $u_\as$ and  $V_\ds$ (both $v_\ds$ and $U_\as$) are s.v. at infinity  (cf. Remark 1.1 of  \cite{Uexit}). 

We present the results of the paper in two subsections below. In the first subsection, we state
our main results  in Theorems  \ref{thm1}  and  \ref{thm2}  and some results complementary to them in Propositions \ref{prop1.1} to \ref{prop1.3}. In the second one, we suppose  $F$ to be attracted to a stable law and present our results as to  asymptotic estimates of  the probability that  $S$ visits $R$ before entering the negative half-line for asymptotically stable walks and obtain asymptotic behaviour of the probability that  $R$ is ever hit by $S$ conditioned to avoid the negative half-line forever.

\v2
\v2

{\bf 1.1.} {\sc Two-sided exit of relatively stable walks.}
\v2

 Let  $a(x)$  be the potential kernel of $S$ when $F$ is recurrent and $G(x)$  the Green kernel when  $F$ is transient :
 $$
 a(x) = \sum_{n=0}^\infty \big[ P[S_n =0] - P[S_n =-x]\big]  \quad\; \mbox{and} \quad\;  G(x) =\sum_{n=0}^\infty P[S_n =x].
  $$
Under the  relative stability $F$ is transient if and only if
$$\int_{x_0}^\infty \frac{H(t)}{A^2(t)}dt <\infty$$
and in this case $A(x)\to\infty$ as $x\to\infty$ (cf. \cite{Urenw}).  


 We shall study  the asymptotic estimate of $P_x[\sigma_R<T]$ or $P_x[\sigma_0 <\sigma_{(R,\infty)}]$. This is not only of  interest in itself but  sometimes useful for the estimates of $P_x(\La_R)$. In fact the comparison of $P_x[\sigma_0 <\sigma_{(R,\infty)}]$ to $1-P_x(\La_R)$ leads to the determination of the asymptotic form  of the renewal sequence---as well as these two  probabilities---under (PRS) with some regularity condition on the negative tail of $F$. 
The following result is fundamental in this direction. Denote by  $g_{\vom}(x,y)$  the Green function of  $S$ killed as it enters $\vom$ (see (\ref{GF}) for the precise definition).

\begin{theorem}\label{thm1} If $F$ is recurrent and p.r.s., then $a(x)$, $x>0$  is s.v. at infinity and as $y\to\infty$ 
$$v_\ds(y) = o\big(a(y)/U_\as(y)\big), \;\; \mbox{and}$$
$$g_{\vom}(x,y) = a(x) - a(x-y) +o\big(a(x)\big)   
 \quad\; \mbox{uniformly for \, $x > y/2$},$$
  in particular $g_\vom(y,y)\sim a(y)$; \;\mbox{and for each constant $\de<1$}
\beqn\label{g/V/a}
g_{\vom}(x,y)  =  
\left\{\begin{array}{ll}
{\displaystyle \frac{V_\ds(x)}{V_\ds(y)} \big[a(y) -a(-y) +o(a(y))\big]} \quad
&\mbox{as\; $y\to\infty$\; uniformly for}\;\; 0\leq x \leq \de y,\\[4mm]
o\big(a(x)U_\as (y)\big/U_\as (x)\big) \;  &\mbox{as\; $x\to\infty$\; uniformly for}\;\;  0\leq y <\de x.
\end{array} \right.
\eeqn
\end{theorem}

 \vskip2mm
 \begin{rem}\label{rem1.1}\, (a)  An intrinsic part of Theorem \ref{thm1} will be proved under a condition weaker than (PRS). (See Proposition \ref{prop3.1} and Remark \ref{rem3.1}(a).) 
 
 (b) Under (PRS) some asymptotic estimates of $a(x)$ and $G(x)$ as $x\to \pm\infty$ are obtained in \cite{Upot} and \cite{Urenw}, respectively; it especially follows that as $x\to\infty$ 
 \beqn\label{1/A}
 1/A(x) \sim \left\{
 \begin{array}{ll}
 a(x)-a(-x) \;\; &\mbox{if $F$ is recurrent},\\[2mm]
 G(x)-G(-x)\quad \; &\mbox{if $F$ is transient}.
 \end{array}\right.
 \eeqn
  (See (\ref{P_a}) and  (\ref{Gpm}) for more about  $a$ and $G$, respectively.)

(c) In view of the identity $g_{\vom}(x,y)=g_{[1, \infty)}(-y,-x) $,  the dual statement of Theorem \ref{thm1} may read  as follows:  If $F$ is recurrent and n.r.s.,
then    as $x\to\infty$  
$$g_{\vom}(x,y) = a(-y)- a(x-y) +o\big(a(-x)\big) \quad \mbox{uniformly for $y>x/2$; \, and}$$
 $$g_{\vom}(x, y)  = 
  \left\{\begin{array}{ll}
  =o\big(a(-y)V_\ds(x)\big/V_\ds(y)\big) \;  &\mbox{as\, $y\to \infty$\;  uniformly for}\;\; 0\leq x<\de y, \\[3mm]
{\displaystyle \frac{U_\as   (y)}{U_\as   (x)} \big[a(-x) -a(x) +o(a(-x))\big]}
&\mbox{as\, $x\to \infty$\; uniformly for}\;\; 0\leq y \leq \de x.
\end{array} \right.
$$

(d)  Condition  (PRS) entails the slow variation  both of $u_\as$ and  of $V_\ds$ (as mentioned previously) which in turn implies  $g_{\vom}(x, y)  = V_\ds (y)U_\as(x)/x\{1+o(1)\} $ for $0\leq x\leq \de y$ (see Remark \ref{rem3.1}(b)). Comparing this to (\ref{g/V/a}) (with $y=x/2$) one sees  that under (PRS)
\[
\begin{array}{ll} 
V_\ds(x)U_\as(x)/x\sim a(x)-a(-x)  \quad\;  &\mbox{if}\quad \limsup a(-x)/a(x)<1,\\[2mm]
V_\ds(x)U_\as(x)/x =o(a(x))  &\mbox{if}\quad \lim a(-x)/a(x)=1.
\end{array}
\]
\end{rem} 
  \v2\v2
  The result corresponding to Theorem \ref{thm1}  for the transient walk is much cheaper---we shall give it in Section 5 as Lemma \ref{lem5.1} in the dual setting (i.e.,  for an n.r.s. walk), whereas  the exact estimation of $u_\as$ (given in Theorem \ref{thm3} below for not $u_\as$ but $v_\ds$ because of the dual setting)  is more costly than for the recurrent walk. 
     Here we state the standard result that if $F$ is transient, then $0<G(0) =1/P[\sigma <\infty] <\infty$ and 
  \beqn\label{gR}
 g_\vom(x,x)  \to G(0).
  \eeqn
(See Appendix (B) for the proof of the latter assertion.)  
  
 If $F$ is n.r.s., then $P_x(\La_R) \to 0$, and  the exact estimation of  $P_x(\La_R)$ seems hard to perform in general.   However, if  the positive and negative  tails are  not balanced in the sense that 
 \beqn\label{e1.7}
  \left\{
\begin{array}{lr}
{\displaystyle  \limsup_{x\to \infty}\,  a(x)/a(-x) <1 }\quad  &\mbox{if \; $F$\, is recurrent,} \\[3mm]
{\displaystyle   \limsup_{x\to \infty}\,  G(x)/G(-x) <1 }\quad\;  &\mbox{if \; $F$\, is transient,} 
    \end{array}\right.
\eeqn
and the positive tail of $F$ satisfies an appropriate regularity condition, then we can compute the precise asymptotic form of $g_\vom(x,y)$ for $0\leq x<\de y$ (that is lacking in the second formula in (c) of Remark \ref{rem1.1}), and thereby obtain that of $P_x(\La_R)$ for  $F$ that is n.r.s. It is pointed out that under (NRS),  this condition entails (\ref{e5.3}) while (\ref{e1.7}) follows from
\[
\left\{
\begin{array}{ll}
{\displaystyle  \limsup_{x\to \infty} \eta_-(x)/\eta_+(x) <1  } \quad  &\mbox{if \; $E|X|<\infty$,} \\[3mm]
{\displaystyle  \limsup_{x\to \infty} A_+(x)/A_-(x) <1  }\quad\;  &\mbox{if \; $F$ is transient,} 
    \end{array}\right. 
\]
where\; $A_\pm(x) =\int_0^x P[\pm X >t] dt$\, ($x\geq 0$) [recall  $\eta_\pm(x) = \int_x^\infty P[\pm X >t] dt$].   (For (\ref{e1.7}), see Appendix (A)---discussed in the dual setting.)We need to  assume 
   \beqn\label{e5.3}
   \exists \lambda>1, \;\;  \limsup \frac{1-F(\lambda t)}{1-F(t)}  <1.
  \eeqn

Put
$$\tilde u(x) = \frac{U_\as(x)[1-F(x)]}{-A(x)}.$$ 
Then if $F$ is n.r.s. and (\ref{e1.7}) holds, 
 \beqn\label{Diff/U}\frac{d}{dt} \, \frac{1}{\hat\ell_\sharp(t)}  =\frac{1-F(t)}{\hat\ell_\sharp^2(t) \hat\ell^*(t)} \sim  \tilde u(t),\eeqn
(see  Remark \ref{rem3.1}(b) and (\ref{dual/sh})) 
so that then
\beqn\label{U'=u}
U_\as(x) \sim \int_0^x \tilde u(t)dt.
\eeqn 

 \begin{proposition}\label{prop1.1}  Suppose that (NRS) holds. Then 
  \beqn\label{LaP1}
P(\La_R) \geq  v_0U_\as(R)[1-F(R)] \{1+o(1)\}. 
 \eeqn
 If (\ref{e1.7}) and  (\ref{e5.3}) hold in addition, then   
  \beqn\label{uPI}
{u_\as(x)} \geq \tilde u(x)\{1+o(1)\},
 \eeqn
 and in case $EX=0$,  
 there exists  a positive constant $c$ such that for sufficiently large  $R$
 $$cP(\La_R) \leq P[\sigma_R <T]  \leq P(\La_R)\{1+o(1)\},$$
 and $u_\as(x) =o\big(U_\as(x)\big/ x\big)$ ($x\to\infty$).
 \end{proposition}
 

In the next theorem we assume, in addition to  (NRS), (\ref{e1.7}) and (\ref{e5.3}),  the continuity of $H_+$:
  \beqn\label{ctn}
 \lim_{\lambda\, \downarrow 1}\limsup_{x\to\infty} \frac{1-F(x)}{1-F(\lambda x)} =1,
  \eeqn
  and  obtain the precise asymptotic forms of $P(\La_R)$ and $u_\as(x)$ in case $EX=0$.
   
   \begin{theorem}\label{thm2} Let $EX=0$ and suppose that  (\ref{ctn}),   (\ref{e1.7}),  (\ref{e5.3}) and (NRS) hold.  Then
   \beqn\label{T2-1}
   u_\as(x) \sim \frac{U_\as(x)[1-F(x)]}{-A(x)}, \quad P(\La_R)/v_0 \sim U_\as(R)[1-F(R)],
   \eeqn
 and for each $\de <1$, uniformly for  $0\leq x< \de y$ as $y\to\infty$
    \beqn\label{T2-2}
    g_\vom(x,y) \sim  \frac{V_\ds (x)U_\as(y)}{|A(y)|(x+1)}\int^{y}_{y-x-1} [1-F(t)]dt; \;\mbox{and}
    \eeqn
    \beqn\label{T2-3}
     \frac{P_{x} [\sigma_y < T]}{P_x(\La_y)} \sim  \frac{a(-y)-a(y)}{a(-y)}.
     \eeqn
\end{theorem}
\v2\v2

\begin{rem}\label{rem5.1.2} (a) \; Under (NRS),  (\ref{e1.7}) follows if we suppose 
\[
\left\{
\begin{array}{ll}
{\displaystyle  \limsup_{x\to \infty} F(-x)]/[1-F(x)] <1  } \quad  &\mbox{if \; $E|X|<\infty$,} \\[3mm]
{\displaystyle  \limsup_{x\to \infty} A_+(x)/A_-(x) <1  }\quad\;  &\mbox{if \; $F$ is transient,} 
    \end{array}\right. 
\]
where $A_\pm(x) =\int_0^x P[\pm X>t] dt$ ($x\geq 0$).  (For (\ref{e1.7}), see Appendix 6.3---discussed in the dual setting.)
\v2
(b) Suppose that $EX=0$,  $\limsup F(-x)/[1-F(x)] <1$ and the \lq continuity' condition (\ref{ctn}) is valid. Then the assumption of Theorem \ref{thm2} holds if
\beqn\label{e1.6}
1-F(x) \asymp L(x)/x \quad (x\to\infty)\quad \mbox{for some s.v.  $L$} 
\eeqn
as is readily verified (the converse is not true). 
\end{rem}

\v2
In case $E|X|=\infty$  one may expect the formulae parallel to those given in Theorem \ref{thm2} to be true
on an ad hoc basis, but  they require  a more delicate analysis than in case $EX=0$. In the next theorem we give a partial result under 
 the following  condition, more restrictive than the continuity condition (\ref{ctn}):
\beqn\label{TRNS}
\exists C>0,\quad p(x) \leq C[1-F(x)]/x \quad\;\;  (x\geq 1).
\eeqn
\begin{theorem}\label{thm3} 
Let $E|X|=\infty$.  Suppose that   (\ref{TRNS}) holds in addition to  (\ref{e1.7}),  (\ref{e5.3}) and (NRS). 
Then the formulae (\ref{T2-1}) and (\ref{T2-2}) hold and  instead of (\ref{T2-3}) it holds that  uniformly for $0\leq x<\de R$,
    \beqn\label{T2-30}
     \frac{P_{x} [\sigma_y < T]}{P_x(\La_y)} \sim  \frac{G(-y)-G(y)}{G(0)}.
     \eeqn
 \end{theorem}
 \v2

 \begin{rem}\label{rem1.2}  Let $E|X|=\infty$. By \cite[Corollary 2]{Ec},  under (NRS) the assumption of  $S$ being oscillating  is equivalent to
 $$\int_1^\infty \frac{1-F(t)}{A_-(t)} dt =\infty.$$
   \end{rem}
 \v2

Because of the identity  $P_{x} [\sigma_y < T]= g_\vom(x,y)/g_\vom(y,y)$, combed with (\ref{T2-2}), Theorem \ref{thm1}, (\ref{gR}) and (\ref{T2-2})   the equivalence (\ref{T2-3})  or (\ref{T2-30})  leads to the following

\begin{corollary}\label{cor1.1} If the  assumption of Theorem \ref{thm2} or that of Theorem \ref{thm3} holds according as $F$ is recurrent or transient, then uniformly for\; $0\leq x<\de R$,
\beqn\label{eC1.1}
P_x(\La_R) \sim |A(R)|g_\vom(x,R) \sim \frac{V_\ds (x)U_\as(R)}{x+1}\int_{R-x-1}^{R}[1- F(t)] dt \qquad \mbox{as $R\to\infty$}.
\eeqn
\end{corollary}

In the dual setting,  Theorem \ref{thm2} is paraphrased as follows. If (PRS) holds,  (\ref{ctn}) and (\ref{e5.3}) hold with $F(-\,\cdot\,)$ in place of $1-F$ and  $\limsup a(-x)/a(x)<1$, then
   $$1-P(\La_y)\sim V_\ds(y)F(-y)\sim v_\ds(y)A(y)   \quad\;  (y\to\infty) ,$$
   $$P_{x} [\sigma_0 <\sigma_{(R,\infty)}] \sim {g_\vom(R,R-x)}/{a(R)} \qquad \mbox{uniformly for } \; 0\leq x\leq R$$
 and for each $\ep>0$,
  \beqn\label{g_prs}
  g_\vom(x,y) \sim  \frac{U_\as (y)V_\ds(x)}{A(x)(y+1)}\int_{x-y-1}^x F(-t)dt \quad (x\to\infty)\;\;\mbox{uniformly for\; $0\leq y< (1-\ep) x$},
  \eeqn
  $$\frac{P_{x} [\sigma_0 <\sigma_{(R,\infty)}]}{P_x[T<\sigma_{(R,\infty)}]} \sim \frac{a(R)-a(-R)}{a(R)}\quad (R\to\infty) \;\; \mbox{uniformly for\; $\ep R < x \leq R $,} 
  $$  
\beqn\label{eC1.10}
\begin{array}{rr}
1-P_x(\La_R) \sim A(R)g_\vom(R,R-x) 
\sim {\displaystyle \frac{U_\as (R-x)V_\ds(R) }{R-x+1}\int_{x-1}^{R} F(-t)dt}\quad  (R\to\infty) 
\\[4mm] 
\quad \mbox{uniformly for\; $\ep R < x \leq R $};
\end{array}
\eeqn
 in particular by (\ref{g_prs}) 
 $$ P_{x} [\sigma_y < T]  \asymp \frac{yV_\ds(x)F(-x)}{V_\ds(y)A(x)}  = o\bigg(\frac{V_\ds(x)/x}{V_\ds(y)/y}\bigg) \quad (x\to\infty) \quad \mbox{uniformly for} \;\; 0\leq y<\de x,$$
 where for the last equality we have used $V_\ds(y)U_\as(y)\asymp y/a(y)$ (see {\bf L(3.1)} in Section 2 and (\ref{A/el/a})).
 It would be plain to state the dual results of Theorem \ref{thm3}. 
\v2
Let $Z$ (resp. $\hat Z$) be  the first ladder height of the strictly ascending (resp. descending ladder) process: $Z= S_{\sigma[1,\infty)}$,  $\hat Z = S_T$. We shall be also concerned with  the overshoot which we define by  
$$Z(R) = S_{\sigma[R+1,\infty)} -R.$$

 \begin{rem}\label{rem1.3}   In \cite[Eq(2.22)]{Uladd} it is shown that if $EZ<\infty$, then
 $$1-P_x(\La_R) \sim \big[V_\ds(R)-V_\ds(x)\big]\big/V_\ds(R) \quad \mbox{ as $R-x \to \infty$ \; for\,  $x\geq0$}, $$
  which, giving an exact asymptotics for $0\leq x \leq\ep R$, partially complements   
 (\ref{eC1.10}). 
  \end{rem}

\v2
From the estimate  of $u_\as$ and $g_\vom(x,y)$  of Theorems \ref{thm2} and \ref{thm3} we can derive some exact asymptotic estimates of 
\beqn\label{LLP}
P_x[ S_{\sigma(R,\infty)-1} =y\,|\, \La_R],
\eeqn
  the conditional probability of  $S$ exiting  the interval $[0,R]$ through  $y$, given $\La_R$. We shall carry out the derivation in Section 6.
Here we state the following consequence of it as to $Z(R)$.
\begin{proposition}\label{prop1.3} Suppose that either the  assumption of  Theorem \ref{thm2} or that of Theorem \ref{thm3} holds. Suppose, in addition,
 that  $1-F(x) \sim L_+(x)/x$ for some s.v. function $L_+$. Then for each $\ep>0$, $\mbox{uniformly for $y>\ep R$ and  $0\leq x <(1-\ep)R$}$, 
as $R\to\infty$
$$ P_x\big[Z(R) \leq y \,\big|\, \La_R\big] \sim 1- \frac{\log [1-(R+y)^{-1}(x+1)]}{ \log[1-R^{-1}(1+x)]}\,;$$
in particular \;
$P_x\big[Z(R) \leq y \,\big|\, \La_R\big] \sim y/(R+y)$\, as\,  $x/R\to 0$.
 \end{proposition}

\v2\v2

{\bf 1.2.} {\sc Asymptotics of $P_x[\sigma_R<T\,|\, \La_R]$ for asymptotically stable walks.}
\v2
 As in \cite{Uexit} we bring in the asymptotic stability condition
\[
{\bf   (AS)} 
\left\{\begin{array}{ll} 
  {\rm  a.} \;\; \;\,  \mbox{$X$ is attracted to  a stable law of exponent $0<\al\leq 2$.} \qquad  \qquad  \qquad  \qquad \\
    {\rm  b.}\mbox{ \; $EX=0$ if $E|X|<\infty$.} \\
  {\rm  c.} \mbox{ \;\;  there exists} \quad 
  \rho :=\lim P[S_n>0].
  \end{array}\right.
  \]

Suppose condition (ASab)---the conjunction of (ASa) and (ASb)---to hold with $\al<2$. It then follows that 
 \beqn\label{stable}
 1-F(x) \sim pL(x)x^{-\alpha}  \quad \mbox{and} \quad F(-x) \sim q L(x)x^{-\alpha}
 \eeqn
 for some s.v. function $L$ and  constant $p=1-q\in [0,1]$.



The next proposition concerns the conditional probability $P_x[\sigma_R <T\,|\,\La_R]$, or what is the same, the ratio $P_x[\sigma_R <T]/P_x(\La_R)$. Note that under (AS), uniformly for $0\leq x\leq R$
\beqn\label{2}
P_x[\sigma_R<T] \sim P_x(\La_R) \sim V_\ds(x)/V_\ds(R) \quad\; \mbox{if} \quad \alpha=2.
\eeqn

\begin{proposition}\label{prop1.4} \, Suppose (AS) to hold with $0<\al<2$ and let $1/2<\de<1$ as above.

{\rm  (i)}\;   For $1< \al <2$, the following equivalences hold:
 \beqn\label{Q_PLa}
   p=0\, \Longleftrightarrow\, (\ref{fun_res}),\, 
 \eeqn
 \beqn\label{eP1.4}   
 pq=0   \Longleftrightarrow \, \lim_{R\to\infty} \frac{P_x[\sigma_R <T\,]}{P_x(\La_R)} =1 \quad\;  \mbox{for some/all\, $x\in \mathbb{Z}$},
 \eeqn
(in the latter case the last limit is uniform for $0\leq x\leq R$), and  if $pq>0$,  then  for some constant  $\theta \in (0,1)$, $P_x[\sigma_R<T] \leq \theta P_x(\La_R) $ $(0\leq x<\de R)$,
and   
  $$ P_x[\sigma_R<T]\sim f\bigg(\frac{x}R\bigg)\frac{V_\ds(x)}{V_\ds(R)} \quad \;\mbox{uniformly for $\; 0\leq  x \leq R$ \; as \; $R\to\infty$}$$
with some increasing and continuous function  $f$ such that  $f(1)=1$ and $f(0) = (\al-1)/\al\hat\rho$. 
     
{\rm  (ii)}\; If $\al=1$, $\rho>0$ and $S$ is recurrent (necessarily $p\leq q$), then
   $$\frac{P_x[\sigma_y<T\,]} {P_x(\La_y)} \,\longrightarrow\,  
    \left\{\begin{array}{ll}
     (q-p)/q \quad &\mbox{as \; $y\to\infty$\; uniformly for}\;\; 0\leq x < \de y,\\[1mm]
0 \;  &\mbox{as \; $x\to\infty$ \;uniformly for}\;\; 0\leq y<  \de x.
\end{array} \right.
$$

 {\rm  (iii)}\; If $\al=1$,  $EX=0$ and $p> q$ (entailing $\rho=0$), then  
 $$\frac{P_x[\sigma_y<T\,]} {P_x(\La_y)} \,\sim\, 
  \left\{\begin{array}{ll}
  (p-q)/p  &\mbox{as \, $y\to\infty$\; uniformly for} \;\;\; 0\leq x \leq \de y,\\[2mm]
  {\displaystyle \frac{p-q}{p}\cdot\frac{U_\as (y)}{U_\as (x)} } \;\; &\mbox{as \; $x\to \infty$\; uniformly  for \; $0\leq y\leq \de x$}.
  \end{array} \right.$$
      
{\rm  (iv)}\;  If  $\rho >0$ and $S$ is transient (necessarily $\alpha \leq 1$), then uniformly for $0\leq x<\de R$,
   $$ P_x[\sigma_R<T\,]/P_x(\La_R) \to 0.$$
 \end{proposition}
 
 %
 
  \vskip2mm
 \begin{rem}\label{rem1.4}   
 (iii) above is obtained as a special case of Theorem \ref{thm2} and does not follow from (ii) by duality. Its proof,  much more involved than that of (ii), crucially depends on the fact that if $p>q$, $P(\La_R)$ is comparable with $P[\sigma_{R}<T]$ and the latter is expressed as $v_0u_\as (R)/g_\vom(R.R)$.   In case $\al=\hat\rho=2p=1$, one may reasonably expect that  $P_x[\sigma_R <T]/P_x(\La_R)$ converges to $0$ (as $R\to\infty$)  whether  $S$ is recurrent or not.
  If $F$ is recurrent (transient),  this were true if we could show that  $P_x[S_{(R,\infty)}> (1+\ep) R\, |\,\La_R] \to 1$  ($P_x[S_{(R,\infty)}> R+1/\ep\, |\,\La_R] \to 1$) as $R\to \infty$, $\ep\downarrow 0$ in this order, which is quite plausible but seems hard to show (even under (AS)).
 \end{rem}


\v2
 
 \v2
 Let  $P^\vom_x$, $x\geq 0$ be   the law of the Markov chain  determined by 
\beqn\label{CndW}
P^\vom_x[S_1= x_1, \ldots, S_n=x_n] =P_x\big[S_1= x_1, \ldots, S_n=x_n,  n >T\big] \frac{V_\ds(x_n)}{V_\ds(x)} 
\eeqn 
  ($x, x_1, \ldots, x_n \geq0$),
  in other words, $P^\vom_x$  is  the $h$-transform with $h=V_\ds$ of the  law of $S$ killed as it enters $\vom$;
  $P^\vom_x$ may be considered to be  the law of  $S$ started at $x\geq0$ conditioned never to enter $\vom$.
From the defining expression (\ref{CndW}) one deduces  that 
 \beqn\label{CndW2}P_x^{\vom}[\sigma_y <\infty] = \frac{V_\ds(y)}{V_\ds(x)}  P_x[\sigma _y <T]
 \quad\;  (x\geq 0, y\geq 0).
 \eeqn
Because of this identity,  the estimates of $P_x[\sigma_R<T]$   obtained in \cite{Uexit} as well as  in this paper lead to the following 
 \v2
 \begin{corollary}\label{cor1.2} Suppose (AS) to hold.  
 \v2
 {\rm (i)}  If  $1<\alpha\leq 2$, then uniformly in $x\geq 0$,
  $$
 P^\vom_x[\sigma_R<\infty] \sim
 f(x/R) 
 $$
where  $f(\xi)$ is a continuous function of $\xi\geq0$ such that
$$
 \begin{array}{ll}
\mbox{for $\xi \leq 1$: } \;  & \left\{\begin{array}{ll}
 \mbox{$f$  is identical to  $1$  if \,  $\alpha=2$,}\\
\mbox{\it $f$ equals  the function appearing in Proposition \ref{prop1.4} if } \; \alpha <2,  
  \end{array}\right. \quad \\[4mm]
\mbox{for $\xi>1$:} \;  &f(\xi) =(\alpha-1) \xi^{-\alpha\hat\rho} \int_0^1 t^{\alpha\rho-1}(\xi -1+t)^{\alpha\hat\rho-1}dt;
 \end{array}
 $$
in particular \;    $P^\vom_x[\sigma_R<\infty]  \sim [(\alpha-1)/\alpha\rho]R/x $\,  as\,  $x/R \to \infty$.
 \v2
 {\rm (ii)}  If  $\alpha=1$ and $F$ is recurrent, then for each $\de<1$,  uniformly in $x\geq 0$ as $R\to\infty$  
 $$
 P^\vom_x[\sigma_R<\infty] 
 \left\{
 \begin{array}{ll}
 \to (q-p)/q   \quad &\mbox{if \;}  q\geq p, \\[2mm]
 \asymp  {\displaystyle \frac{R(1-F(R)]}{-A(R)} }\to 0 \quad &\mbox{if \;}  \, q < p, 
 \end{array}\right. \;\;\;  \mbox{for}\;\;  x <\de R,
 $$
 $$
 P^\vom_x[\sigma_R<\infty] 
 \left\{
 \begin{array}{ll}
 \to 0   \quad &\mbox{if \;}  \, q\geq p, \\[2mm]
 \sim {\displaystyle \frac{p-q}{p}\cdot \frac{R/A(R)}{x/A(x)} } \quad &\mbox{if \;}  \, q < p. 
 \end{array}\right.\;\;\; \mbox{for}\;  x> R/\de.
 $$
 \v2
 {\rm (iii)}  If $F$ is transient, then for each $\de<1$, as $R\to\infty$
$$ 
 P^\vom_x[\sigma_R<\infty] \to 0  \quad \mbox{uniformly for\; $0\leq x<\de R$\;  and as} \;\;   x-R \to\infty.$$
 \end{corollary}
 \v2\noindent
 \pf In view of (\ref{CndW2}),  (i) follows from Proposition \ref{prop1.4}(i) in case  $x\leq R$ and from Lemma \ref{lem7.3} in case  $x>R$. As for (ii),  use (ii) and (iii)  of  Proposition \ref{prop1.4} together with the estimate of $P_x(\La_R)$ of Proposition \ref{prop1.1} (in case $q<p, x<\de R$).   (iii) follows from  Proposition \ref{prop1.4}(iv) in case  $\rho>0$ and from 
{\bf L(4.5)} of the next section in case $\rho=0$. \qed  
\v2

\begin{rem}\label{rem1.5}   If $\sigma^2<\infty$, the relations given in (\ref{2})  and Corollary \ref{cor1.2} for $\alpha =2$ are valid.  The asymptotic form of $P_x^{\vom}[\sigma_R<\infty]$ for $\alpha=2$ as $x/R\to\infty$ follows from the invariance principle for a random walk  conditioned to stay positive as established in \cite{CC}, but the validity of the corresponding statement is not clear for the case  $1<\alpha<2$.
\end{rem}

The rest of this paper is organised as follows. In Section 2, we state some of the results from \cite{Uexit} and some known facts that are fundamental in the later discussions. Proof of Theorem \ref{thm1} is given in Section 3. In Section 4 we prove Proposition \ref{prop1.1} in case $EX=0$ and Theorem \ref{thm2} after showing miscellaneous lemmas  in preparation for the proofs. Proposition \ref{prop1.1} (in case $E|X|=\infty$) and Theorem \ref{thm3} are proved in Section 5. In Section 6, we derive asymptotic estimates of the conditional probability in (\ref{LLP}). In Section 7 we deal with asymptotically stable walks and prove Proposition \ref{prop3.1};
 for the proof we compute, in Lemma \ref{lem*1}, the exact asymptotic forms of the renewal sequences $v_\ds$ and $u_\as$ that are of independent interests. 



\section{Preliminaries}

By the fact that $V_\ds$ is harmonic for the r.w. killed as it enters $\vom$  we have
 \beqn\label{upp_bd}
P_x(\La_R) \leq V_\ds(x)/V_\ds(R).
\eeqn 
(see \cite[Eq(2.3)]{Uexit}).
If either  $V_\ds$ or $\int_0^x P[\hat Z>x]$ is regularly varying,  it follows \cite[Eq(8.6.6)]{BGT}  that
 \beqn\label{U/Z}
 \frac{V_\ds(x)}{x v_0}\int_0^xP[-\hat Z>t] dt\, \longrightarrow\, \frac1{\Ga(1+\alpha\hat\rho)\Ga(2-\alpha\hat\rho)}.
 \eeqn

For a non-empty $B\subset \mathbb{Z}$ we define  the Green function  $g_B(x,y)$ of the r.w. killed as it hits $B$ by
\beqn\label{GF}
g_B(x,y) = \sum_{n=0}^\infty P_x[ S_n=y, n<\sigma_B].
\eeqn
(Thus if $x\in B$,  $g_B(x,y)$ is equal to $\de_{x,y}$ for $y\in B$ and to $ E_x[g_B(S_1,y)]$ for $y\notin B$.)  We shall repeatedly apply  the formula
\beqn\label{G-ft}
g_{\vom}(x,y) =\sum_{k=0}^{x\wedge y} v_\ds(x-k)u_\as   (y-k)\quad \mbox{for} \; x,y \geq 0
\eeqn
\cite[Propositions 18.7,  19.3]{S}.  



We shall use several results from \cite{Uexit}. 
Here are given some of those that are of the repeated use.


{\bf L(2.1)}\,   For $0\leq x\leq R$, \; 
${\displaystyle \sum_{y=0}^R g_{\vom}(x,y) \leq  V_\ds(x)U_\as(R).}$
\v2\v2\noindent
  We have shown (\ref{fun_res}) under the following condition (among others):
\v2
(C3) \quad both $V_\ds(x)$ and $xU_\as(x)$  are s.v. as $x\to\infty$.
\v2\noindent
We shall need the dual results of those valid under (C3) whose dual is give as:
\v2
$(\widehat{{\rm C}3})$ \quad both $xV_\ds(x)$ and $U_\as(x)$  are s.v. as $x\to\infty$.
\v2\noindent
The condition  (C3) follows from (PRS) and $(\widehat{{\rm C}3})$ from (NRS) as mentioned previously.
 Put for $t\geq0$
\beqn\label{3.0}
\ell^*(t) = \int_0^tP[Z>s]ds, \quad\mbox{and}\quad  \hat\ell^*(t) = \frac1{v_0}\int_0^tP[-\hat Z>s]ds
\eeqn
(as in \cite{Uexit})  and 
\beqn\label{3.01}
\ell_\sharp(t) = \int_t^\infty \frac{F(-s)}{\ell^* (s)}ds \quad\mbox{and} \quad \hat\ell_\sharp(t) = \int_t^\infty \frac{1-F(s)}{\hat\ell^* (s)}ds 
\eeqn
(slightly differently from \cite{Uexit} in case (C3) or $(\widehat{{\rm C}3})$ fails: see (\ref{el/sh})). It is known that  $Z$ is r.s. if and only if $xU_\as(x)$ is s.v. which in turn is equivalent to the slow variation of $\ell^*$ is s.v. \cite{R}.
\cite[Appendix(B)]{Upot}, \cite{Urenw}.
\v2
{\bf L(3.1)}  Under (C3), $\ell^*$ and $\ell_\sharp$ are s.v.,\; $u_\as(x) \sim 1/{\ell^*(x)}\quad 
\mbox{and} \quad V_\ds(x)\sim 1/\ell_\sharp(x).$
\v2\noindent
By the duality this entails that under $(\widehat{{\rm C}3})$,\, $\hat\ell^*$ and $\hat\ell_\sharp$ are s.v.,
\beqn\label{dual/sh}
v_\ds(x) \sim 1/\hat\ell^*(x) \quad\mbox{and} \quad U_\as(x) \sim 1/\hat\ell_\sharp(x). 
\eeqn
\v2
{\bf L(3.3)}  \;If either (C3) or $(\widehat{{\rm C}3})$ hold, then
\;  $V_\ds(x)U_\as   (x) H(x) \,\longrightarrow \, 0.$
\v2
\v2

 
 {\bf L(3.4})\,  If  (C3) holds,  then for each  $\ep >0$. $P_x[Z(R) > \ep R\,|\, \La_R] \to 0$\; ($R\to\infty$)\; uniformly for $0\leq x<R$.


\v2\v2


{\bf L(4.5)}\;  If either $(\widehat{{\rm C}3})$ or (AS) with $\alpha<1=\hat \rho$ holds, then 
$P_x(\La_R)V_\ds(R)/V_\ds(x) \,\longrightarrow \, 0 \;\; (R\to\infty) \;\;  \mbox{uniformly for $0\leq x <\de R$.}
$
\v2
These results follow from Lemmas 2.1, 3.1, 3.3, 3.4 and 4.5 of \cite{Uexit}.

\section{Proof of Theorem \ref{thm1} and related results.}  




\v2
Here we shall suppose that $F$ is recurrent.
For the present purpose it is convenient to consider  the Green function of $S$ killed as it hits $(-\infty,0]$,  instead of   $\vom =(-\infty,-1]$.
 We make this choice for convenience in applying the identity (\ref{ladd})  below. Put
$$g(x,y)  = a(x) +a(-y)-a(x-y).$$
Then, for $x\geq 1$, 
\beqn\label{ladd}
E_x[a(S_{\sigma_{(-\infty,0]}})] = a(x) - V_\ds(x-1)/EZ, \; \mbox{and} 
\eeqn
\beqn\label{g/g}
g_{(-\infty,0]}(x,y) = g(x,y) - E_x[g(S_{(-\infty,0]}, y)],
\eeqn
which take less simple  forms  for  $E_x[a(S_{T})] $ and $g_{\vom}(x,y)$.
Here (\ref{ladd}) follows from Corollary 1 of \cite{Uladd} and (\ref{g/g}) from the identity
 $g_{\{0\}}(x,y) = g(x,y)$ ($x\neq 0$) (cf. \cite[P29.4]{S}). 


  
  We bring in the following conditions:
   \beqn\label{Ca-c}
   \begin{array}{ll}
  {\bf   (1)} \quad   \mbox{$a(x)$  is  almost increasing and  $a(-x)/a(x)$ is bounded as $x\to\infty;$} \qquad\qquad\\[1mm]
 {\bf   (2)} \quad 
 {\displaystyle \sup_{x: -z \leq x\leq \de z} \frac{|a(x-z)-a(-z)|}{a(z)} \, \longrightarrow \, 0} 
\quad \mbox{as $z \to \infty$ \;    for any $\de <1;$} \\[4mm]
 {\bf   (3)} \quad P_x[S_T> -\ep x] \to 0 \quad \mbox{as\;  $x\to\infty$\, and\, $\ep\downarrow 0$\; in this order.}
 \end{array}
  \eeqn
  These  are all satisfied if either (PRS) or (AS) with  $\alpha=1$ and $\rho>0$ holds (see Remark \ref{rem3.1}(a),(c) below). 
  \begin{proposition}\label{prop3.1} 
 Suppose conditions  (1) to (3) above to hold. Then it holds that
  $$\limsup_{x\to\infty} a(-x)/a(x)\leq 1$$
   and that
for any $\ep>0$,  as\; $x\to\infty$
 \[
   \begin{array}{ll}
{\rm (i)}\;\;\; g_{\vom}(x,y)=
  a(x) -a(x-y) + o(a(x)) \quad  \mbox{ uniformly for}\;\; x> \ep y > 0; \\[2mm]
  {\rm (ii)} \;\; g_{\vom}(x,y)= a(x) -a(-y) +o(a(y)) \quad  \mbox{uniformly for}\;\; \ep y\leq x \leq (1-\ep) y;
\end{array}
\]
 in particular\,  $g_\vom(x,y) \sim a(x)$  if $a(x-y) =o(a(x))$ and $x>y/2 > 0$. 
\end{proposition}
  \v2\v2
  
  \begin{rem}\label{rem3.1}    (a)  Let (PRS) hold.  Then according to \cite[Theorem 7]{Upot},  it hods  that 
  \beqn\label{P_a}
  \begin{array}{ll} 
  a(x)\;\mbox{is s.v.}; \quad a(x)-a(-x) \sim 1/A(x); \\[2mm]
  {\displaystyle a(x) \sim \int_0^x \frac{F(-t)}{A^2(t)}dt, \quad  a(-x) = \int_0^x \frac{1-F(t)}{A^2(t)}dt +o(a(x))}, \qquad
  \end{array}
  \eeqn
as $ x\to\infty$. Combined with (PRS)  these  yield  that for $- z\leq  x< 0$
  $$a(x-z) -a(-z) = \int_{-z}^{x-z} \frac{F(-t)}{A^2(t)}dt +o(a(z)) = o\bigg(\frac{x}{A(z)z}\bigg) + o(a(z)),$$
and similarly for  $0<x< z$, and one sees that  (2) of (\ref{Ca-c}) is satisfied.  We also know that  by  {\bf L(3.1)} $V_\ds(x)$ is s.v.  so that for each $M>1$, $\lim P_x[S_T < -M x] =0$, in particular (3) of (\ref{Ca-c}) is satisfied. (1) is obvious  from (\ref{P_a}).

(b) Let (PRS) hold. By virtue of the Spitzer's formula (\ref{G-ft}) we know $g_\vom(x,y)\sim U_\as (x)/\ell^*(y)$ for $0\leq x< \de y$ ($\de<1$), which combined with (\ref{P_a}), $V_\ds(x)\sim 1/\ell_\sharp(x)$ and the second half of Proposition \ref{prop3.1} leads to the equivalence relations
\beqn\label{aaa}
a(x)\asymp \frac1{A(x)} \, \Longleftrightarrow \, \limsup\frac{a(-x)}{a(x)}<1  \, \Longleftrightarrow \, a(x)\asymp \frac1{\ell^*(x)\ell_\sharp(x)},
\eeqn 
as well as  $1/\ell^*(x)\ell_\sharp(x) =g_\vom(x,x/\de)= a(x)-a(-x) +o(a(x))$, so that 
each of the conditions in (\ref{aaa}) implies
\beqn\label{A/el/a}
A(x) \sim \ell^*(x)\ell_\sharp(x); 
\eeqn
and  if  $\lim a(-x)/a(x)=1$, then  $\ell^*(x)\ell_\sharp(x)a(x) \to \infty$. In  Appendix (A) we shall present a sufficient condition expressed in terms of integrals of $F$ for (\ref{A/el/a})  to hold under (C3). 

(c)\, Suppose that $F$ satisfies (AS) with  $\alpha=1$ and $\rho> 0$. Then conditions (1) to (3) are satisfied.
 Indeed, for $\rho=1$, (PRS) is satisfied,  while for $0<\rho<1$, 
$$a(x)\sim a(-x) \sim c_\rho  \int_1^x \frac{dt}{tL(t)}$$
 with a certain positive constant $c_\rho$, according to \cite[Proposition 61(iv)]{Upot}, entailing (1) and (2) of (\ref{Ca-c}); moreover $V_\ds$ is regularly varying  with index $1-\rho\in (0,1)$ so that (3) of (\ref{Ca-c}) is satisfied owing to the generalised arcsin law \cite[p.374]{F}. Noting (\ref{P_a}) is valid if $\rho=1$, in  view of Proposition  \ref{prop3.1}  these also show that for each $\ep >0$,
  \[g_{\vom}(x,y)= a(y) -a(-y) +o(a(y))\quad \; \mbox{uniformly for}\;\; \ep y\leq x \leq (1-\ep) y.
\]
 Below we verify that
\beqn\label{g/SC}
g_\vom(x,y) =o(a(y)) \quad (y\to\infty) \quad \mbox{uniformly for} \;\; x>(1+\ep)y.
\eeqn
 To this end  
we  have only to show $\sup_{x>(1+\ep)y}P_x[\sigma_y<T] \to 0$. This is immediate,  if $\rho=1$,  for then $P_x[S_T<-Mx]\to 1$ as $x\to \infty$ for each $M>1$.
For $0<\rho<1$, by (3)  one observes that
  $$P_x[\sigma_y < T] =  P_x[\ep_1<S_{\sigma[0, y]}/y<1-\ep_1 ] 
  \sup_{z: \ep_1 <z/y<1-\ep_1} P_z[\sigma_y< T]+ o_{\ep_1}(1)$$
with $o_{\ep_1}(1) \to 0$ as $\ep_1\to 0$, but $P_z[\sigma_y< T]= g_{\vom}(z,y)/g_{\vom}(y,y)$ tends to zero uniformly in $z$ since  $a(y)\sim a(-y)$ and  $\bar a$ is s.v.  Thus (\ref{g/SC}) is verified.
We shall derive in Section 6  essentially the same result as in Proposition \ref{prop3.1} but under (AS)  in the case $\alpha =1$ and  $0<\rho<1$  except for $\rho=1/2$, so that the inclusion of that case 
is significant.   \end{rem}

\v2

We state the following corollary that follows immediately from Proposition \ref{prop3.1} and Remark \ref{rem3.1}(c) because of
\beqn\label{Id}
P_x[\sigma_y < T] = \frac{g_\vom(x,y)}{g_\vom(y,y)}
\eeqn
as well as the fact that if $V_\ds$ is s.v., then $P_x[\sigma_{(-\infty,-R]} =T] \to 1$ as $x/R \to\infty$.
 \begin{corollary}\label{cor3.1}
 Suppose that  $F$ is recurrent and satisfies either (PRS) or (AS) with $\alpha=1$, $\rho>0$. Then,  for any $0<\de<1$, as $R\to\infty$ 
$$ P_x[\sigma_R<T\,] = \left\{ 
\begin{array}{ll}
{\displaystyle \frac{a(R)-a(-R)}{a(R)} +o(1) }\quad &\mbox{uniformly for \; $\ep R\leq x< \de R$},\\
o(1)  &\mbox{uniformly for \; $x> R/\de$}.
\end{array}\right. $$

[For a transient r.w. $a(x)= G(0)- G(-x)$ so that $a(x)-a(-x)=G(x)-G(-x)$.  If  a transient $F$  is p.r.s., by the estimate of $G(x)$  given in  \cite{Urenw} (see (\ref{Gpm}),  Lemma \ref{lem5.1}---given in dual setting---of this chapter)  one sees that a formula analogous to the above holds but for  $R/\de <x<R/\ep$ (resp. $x< \de R$) instead of $\ep R\leq  x<\de R$ (resp. $x>R/\de$) .]
\end{corollary}
\v2

The proof of Proposition \ref{prop3.1} and Theorem \ref{thm1}  will be given after showing two lemmas.

\begin{lemma}\label{lem3.1}\, If (1) and (2) of (\ref{Ca-c}) hold,  then for any $0<\de<1$ and $M\geq1$,
\beqn\label{L4.1}
\sup\bigg\{\frac{|a(x-z) -a(-z)|}{a^\dagger(|x|) +a(-z)} : - Mz <x< \de z\bigg\} \,\longrightarrow\, 0 \qquad (z\to\infty),
\eeqn
where $a^\dagger(0) =1$ and $a^\dagger(x) = a(x)$ if $x\neq 0$.
\end{lemma}
\pf  From (1) and (2) it follows that 
\beqn\label{4.6}
|a(x-z) -a(-z)|/a(|x|) \to 0\quad (z\to\infty) \quad\mbox{uniformly for $-Mz<x< -z$}
\eeqn  
 Indeed, if $-2z \leq x < - z$,  putting $x'= - z$, $z' = x'-z=-2z$ and writing
$$a(x-z)- a(-z) = [ a(x-z')-a(z')] + [a(x'-z)-a(-z)]$$
one sees, using the sub-additivity of $a(\cdot)$,  that $|a(x-z)- a(-z)|/a(|x|) \to 0$ under (1) and (2); for the general case divide the interval $[x-z]$ at  multiples of $-kz$ and  write  $a(x-z)- a(-z)$ as a telescopic sum. 
 
Pick $\ep>0$ arbitrarily and choose $z_0$---possible under (2)---so that 
\beqn\label{aza}
|a(x-z) -a(-z)| < \ep^2 a(z) \quad \mbox{ whenever  $z\geq z_0$ and $-z\leq x \leq \de z$}.
\eeqn
  Let $-z\leq x<0$ and $z>z_0$. On the one hand, if $a(-z) \leq \ep a(z)$, by the inequalities
$$-\frac{a(-z)}{a(z)} a(x) \leq a(x-z)-a(-z) \leq \frac{a(x-z)}{a(z-x)} a(-x)$$
(cf. \cite[Lemma 3.2]{Uladd},\cite[Section 7.1]{Upot}), condition (1) entails that $|a(x-z)-a(-z)| \leq \ep Ca(|x|)$, where we have also used the bound  (\ref{aza}) to have $a(x-z)/a(-x+z) \leq C' [a(-z) +\ep^2a(z)]/ a(z)   <2C'\ep.$  On the other hand if  $a(-z) > \ep a(z)$, (2) entails $|a(x-z)-a(-z)| <\ep^2 a(z)< \ep a(-z)$. Thus  we have $|a(x-z)-a(-z)| \leq \ep C[a(|x|) + a(-z)]$, showing  the supremum in (\ref{L4.1}) restricted to  $-z \leq x\leq 0$ tends to zero.  

The same argument as above applies to the case $0<x<\de z$ 
 by noting  that $a(z-x)\geq > c a(z)$ with $c>0$ and   $a(x-z) = a(-z)   +o\big(a(z)\big)$
because of the sub-additivity  and  of  (2), respectively.
\qed

 \begin{lemma}\label{lem3.2} If (1) to (3) of (\ref{Ca-c}) hold, then for any $\ep>0$, as  $y\to\infty$
 \beqn\label{H-y}
 E_x[a(S_{T}-y)]  =  E_x[a(S_{T})]  + o(a(x))\quad \mbox{uniformly for}\; \; x>\ep y.
 \eeqn
 \end{lemma}
 \pf  By (3), for any $\ep>0$ and $\ep_1>0$ we can  choose $\de>0$ so that 
 \beqn\label{*}
 P_x[S_T \geq -\de y] <\ep_1 \quad \mbox{if}\;\; x>\ep y.  
 \eeqn
 Supposing (1) and (2) to hold we  apply Lemma \ref{lem3.1} to see that as $y\to\infty$,  $$a(-z-y) = a(-z) +o\big(a(-z)+a(y)\big) \quad \mbox{uniformly  for}\; \; z>\de y,$$
 hence
 \beq
 E_x[a(S_T-y); S_T <-\de y] &=& E_x[a(S_T); S_T < -\de y]\{1+o(1)\} + o(a(y))\\
 &=& E_x[a(S_T)] \{1+o(1)\} + O\big(\ep_1 a(y)\big) +o(a(y)).
 \eeq
 Here (\ref{*}) as well as (1) is used for the second equality, and  by  the same reasoning (with the help of (2)) the left-most member is  written as $E_x[a(S_T-y)] + O(\ep_1 a(y))$.  Noting   $ E_x[a(S_T)] \sim E_x[a(S_{\sigma_{(-\infty,0]}})] \leq C a^\dagger (x)$, we can conclude that  
 $E_x[a(S_T-y)]= E_x[a(S_T)]+o(a(x\vee y))+  O(\ep_1 a(y))$, which shows (\ref{H-y}),
 $\ep_1$ being arbitrary and $a(\cdot)$ sub-additive. \qed
 \v2
 
 {\bf Proof of Proposition \ref{prop3.1}.}  Note that for the asymptotic estimates under $P_x$, $T$ and $\sigma_{(-\infty,0]}$ may be interchangeable as $x\to \infty$. Then 
 applying  (\ref{g/g}) and Lemma \ref{lem3.2} in turn one sees that as $y\to\infty$
  \beq
  g_{(-\infty,0]}(x,y) &=& a(x) - a(x-y) - \big( E_x[a(S_{(-\infty,0]})] - E_x[a(S_{(-\infty,0]}-y)] \big) \\
  &=& a(x) -a(x-y) +o(a(x))
  \eeq
  uniformly for $ x>\ep y$, showing  the fist formula of the proposition.  Since   (2) of  (\ref{Ca-c}) entails $a(x-y) - a(-y) =o(a(y))$ for $x<(1-\ep) y$,
  the second formula  follows. Combining (i) and (2) one easily verifies $\limsup a(-x)/a(x)\leq 1$.
  \qed
 
  \v2

 {\bf Proof of Theorem \ref{thm1}.}  Taking (\ref{fun_res})
 into account we show that under (PRS)  for any $0<\de<1$
 $$g_{\vom}(x,y)  = 
  \left\{\begin{array}{ll} P_x(\La_y)\big[a(y) -a(-y) +o(a(y))\big] \quad
&\mbox{uniformly for}\;\; 0\leq x  < \de y,\\[1mm]
o(a(y)) \;  &\mbox{uniformly for}\;\; x > y/\de.
\end{array} \right.
$$
We have only to consider the case $x = o(y)$, since  the other cases  readily follows from Proposition \ref{prop3.1} because of the slow variation of $V_\ds$ and  of $a$ (see (\ref{P_a}) for the latter). However,  uniformly for $0\leq x < \frac14 y$ we have
\beq
g_{\vom}(x,y) &=& \sum_{\frac14 y\leq w<\frac12 y} P_x[S_{\sigma[\frac14 y,\infty)}=
w]g_\vom(w,y) + o\big(P_x(\La_{\frac14 y})g_\vom(y,y)\big)\\
 &=& P_x(\La_{y})[a(y)-a(-y)+o(a(y)],  
\eeq
where the first and second equalities are due  to {\bf L(3.4)} and Proposition \ref{prop3.1}, respectively. 

With the help of  $g_{-\vom}(-x,-x) \sim a(x)$ (valid under (PRS)) and $g_{\vom}(x, 0) = v_\ds   (x)$ ($x\geq 1$), one can easily deduce from  {\bf L(4.5)} that
\beqn\label{v/aU}
v_{\ds}(x) =o\big(a(x)/U_\as (x)\big).
\eeqn
(the dual assertion, given by (\ref{u/aV}), is more naturally derived). 
The second formula of (\ref{g/V/a}) follows from (\ref{v/aU}) in view of the Spitzer's formula (\ref{G-ft}) for $g_\vom(x,y)$.  \qed 
  
For later usage, here we state the  result corresponding to Corollary \ref{cor3.1} for the n.r.s. walk (see Remark \ref{rem1.1}(c)):, then 
\beqn\label{Dual} P_x[\sigma_y<T\,] = \left\{ 
\begin{array}{ll}
o\big(V_\ds(x)\big/V_\ds(y) \big)  &\mbox{as\; $y\to\infty$\; uniformly for\; $0\leq x< \de y$},\\[1mm]
  {\displaystyle \frac{U_\as (y)}{U_\as (x)}\bigg[\frac{a(-x)-a(x)}{a(-x)} +o(1) \bigg] } \;\; 
  &\mbox{as\; $x\to\infty$\;  uniformly for }\;    y <\de x.\\
\end{array}\right. 
\eeqn


\section{Proof of  Proposition \ref{prop1.1} (case $E|X|<\infty$)\\ and Theorem \ref{thm2} }

This section consists of three subsections.  In the first one we obtain some basic estimates of $P(\La_R)$ and $u_\as (x)$ for negatively relatively stable walks and prove Proposition \ref{prop1.1}. In the second we verify precise asymptotic forms of $P(\La_R)$ and $u_\as (x)$ asserted in  Theorem \ref{thm2}. The rest of  Theorem \ref{thm2}  is proved in the third one.  
We use the following notation: 
$$H_+(x) :=1-F(x)= P[X > x] \quad (x\geq0);$$
$$B(R) := \mathbb{Z} \setminus [0,R];$$ 
    $$N(R) =\sigma_{B(R)} -1,$$
    the first time the r.v. leaves $B(R)$ after time zero. Throughout this section we suppose  
$$xH(x)/A(x)\to -\infty  \qquad (x\to\infty).
$$
Recall that this entails $\rho=0$, $V_\ds(x)\sim x/\hat\ell^*(x)$ and $U_\as(x)\sim 1/\hat\ell_\sharp(x)$. 
 \v2\v2
 {\bf 4.1.} {\sc  Preliminary estimates of $u_\as $ and  proof of Proposition \ref{prop1.1} in case $EX=0$.}
 \v2

In this subsection, we are mainly concerned with  the case  $F$ being recurrent but some of the results are valid also for $F$ being transient.  Note that the potential function $a(x)$ is always well defined and approaches the constant $G(0)$ as $|x|\to\infty$ if $F$ is transient.

Since $g_\vom(0,y)= v_0u_\as   (y)$,  Theorem  \ref{thm1} shows that  if (NRS) holds, then   
\beqn\label{sigm/ua}
P[\sigma_y< T] \sim g_\vom(0,y)/ a(-y) = v_0 u_\as (y)/a(-y),
\eeqn
 hence by the dual of  {\bf L(4.5)}
\beqn\label{u/aV}
u_\as (y) =o\big(a(-y)\big/V_\ds(y)\big).
\eeqn
\v2

\begin{lemma}\label{lem4.1} Suppose that (NRS) holds and $\limsup a(x)/a(-x)<1$ (necessarily  $EX=0$). Then
\beqn\label{u/U}
u_\as   (y) = o\Big(1\big/\big[V_\ds(y)|A(y)|\big]\Big).
\eeqn
\end{lemma}
\pf  The second condition of the supposition implies $a(-y) \asymp 1/A(y)$ owing to (\ref{P_a}). Hence (\ref{u/U}) is immediate from (\ref{u/aV}). \qed

 \v2
Under $\limsup a(x)/a(-x)<1$,  by  (\ref{A/el/a})  in Remark \ref{rem3.1}(b)   we have $1/A(x)\sim V_\ds(x)U_\as (x)/x$, so that (\ref{u/U}) is equivalently stated as
\beqn\label{u/U/y}
u_{\as}(y) \sim o(U_\as (y)/y),
\eeqn
 which is an expected consequence, for  $U_\as $ is s.v.;  we shall use  this expression instead of  (\ref{u/U}).




\begin{lemma}\label{lem4.2}   Suppose  (NRS) to hold.  
\v2
{\rm  (i)} \quad $P(\La_R) \geq U_\as    (R)H_+(R)\{v_0+o(1)\}.$

\v2
{\rm  (ii)} \quad If $1-F$ is of dominated variation, then 
for each  $0<\ep< 1$, \; 
$$P\big[\ep R \leq  S_{N(R)} <(1-\ep)R\,\big|\, \La_R\big] \to 0\; \;\;\mbox{and}\;\;\;
P\big[ S_{N(R)} < \ep R, \La_R ] \asymp U_\as (R)H_+(R).$$
[Here a non-increasing function $f$ is of dominated variation if $\liminf f(\frac12 x)/f(x)>0$ \cite{F},\cite{BGT}.]
  \end{lemma}

 \v2
 
\pf 
  We need to find an appropriate upper bound of  $g_{\mathbb{Z}\setminus [0,R)}(x,y)$. To this end
  we use the identity 
\beqn\label{idn_g}
 g_{B(R)}(x,y)  =g_{\vom}(x,y) -  E_x[g_{\vom}(S_{\sigma(R,\infty)},y);\La_R] \quad (0\leq x, y <R).
 \eeqn  
By Spitzer's representation (\ref{G-ft})  we see  that  for  $z\geq R$, $1\leq y<\de R$ ($\de<1$),  $g_{\vom}(z,y) \leq U_\as   (y)/\hat\ell^*(R)\{1+o(1)\}$, so that for  $0\leq x, y<\de R$,
\beqn\label{LB_g}
g_{B(R)}(x,y)  \geq g_{\vom}(x,y) - U_\as   (y)P_x(\La_R)/\hat\ell^*(R)\{1+o(1)\}. 
\eeqn
Since 
$\sum_{y=0}^{R/2} U_\as(y) 
\sim  \frac12\,RU_\as (R)$
and   $R/\hat \ell^*(R) \sim  V_\ds(R)$ and since $P(\La_R)V_\ds(R)\to 0$ by virtue of  {\bf L(4.5)},
using $g_{\vom}(0,y) =v_0 u_\as   (y)$ we accordingly deduce from (\ref{LB_g}) 
  \beqn\label{sum/2}
  \sum_{y=0}^{R/2} g_{B(R)}(0,y) =  U_\as   (R)\{v_0 +o(1)\}.
  \eeqn
Hence  
 \beqn\label{La/U/F}
  P(\La_R) \geq \sum_{y=0}^{R/2} g_{B(R)}(0,y)H_+(R-y)\geq U_\as   (R)H_+(R)\{v_0+o(1)\},
  \eeqn
 showing (i).  
 
 The first probability in (ii) is less than $\sum_{y=\ep R}^{(1-\ep)R} u_\as (y)H_+(R-y)$. After summing by parts, this sum may be expressed as 
 \beqn\label{Diff}
 \Big[U_\as (y)H_+(R-y)\Big]_{y=\ep R}^{(1-\ep) R} -  \int_{y=\ep R}^{(1-\ep)R} U_\as (t)d H_+(R-y)
 \eeqn
apart from the error term of smaller order of magnitude than  $U_\as (R)H_+(\ep R)$.  Because of the slow variation of  $U_\as $ the above difference  is  $o\big(U(R) H_+(\ep R)\big)$. By (i) we  therefore obtain  the first relation of (ii), provided that $1-F$ is of the dominated variation.  The second relation of (ii) follows from the first and (i).  \qed

\begin{lemma}\label{lem4.3} 
 Suppose $\limsup [1-F(\lambda x)]/[1-F(x)] <1$ for some $\lambda >1$. Then \beqn\label{L4.8} P_x[ Z(R) >MR\,|\, \La_R] \to 0 \quad\mbox{as\; $M\to\infty $\; uniformly for}\;\;  0\leq x\leq R. \eeqn 
  \end{lemma}   
  \pf
For any integer  $M>1$, writing $M'=M+1$ we have
$$P_x[Z(R) >  MR\,|\, \La_R] \leq \frac{\sum_{y=0}^{R} g_{B(R)}(x,y)H_+(M'R-y)}{\sum_{y=0}^{R} g_{B(R)}(x,y)H_+(R-y)}\leq \frac{H_+(MR)}{H_+(R)},
$$
of which the last member approaches zero as $M\to\infty$ uniformly in $R$ under the supposition of the lemma.
\qed

\v2
  
 Lemma \ref{lem4.2}(ii) says that, given $\La_R$, the conditional law of $S_{\sigma(R,\infty)-1}/R$, the position of departure of the scaled r.w. $S_n/R$ from the interval $[0,1]$, tends to concentrate  near the boundary.  If the positive tail of $F$ satisfies the continuity condition (\ref{ctn}),  we shall see that such a concentration should be expected to occur  only  about the lower boundary (see Lemma \ref{lem4.5}), otherwise, this may be not true.

 \v2\v2
{\bf Proof of Proposition \ref{prop1.1} (case  $EX=0$).}  Suppose the assumption of Proposition \ref{prop1.1} to hold, namely
\beqn\label{e5.1}
\left\{
\begin{array}{ll}
{\rm (a)}\;\; \mbox{$F$\;  is n.r.s.};\\[1mm] 
{\rm (b)}\;\; EX=0\quad \mbox{and}\quad  {\displaystyle  \limsup_{x\to \infty} {a(x)}/{a(-x)} <1};
\qquad\qquad\qquad\qquad\qquad
\\[1mm]  
{\rm (c)}\;\; \exists \lambda>1, \;\;  \limsup {H_+(\lambda t)}/{H_+(t)}  <1.
\end{array}\right.
\eeqn  
\v2\noindent
 Then Lemmas \ref{lem4.2}(i)  and \ref{lem4.3} are applicable, of which the former one gives the lower bound  of $P(\La_R)$ asserted in Proposition \ref{prop1.1}. 
\


 By Lemma \ref{lem4.2}(i)   we have 
  $$P(\La_R)/v_0  \geq U_\as(R)H_+(R)\{1+o(1)\},$$
   the lower bound  of $P(\La_R)$ asserted in Proposition \ref{prop1.1}.
By (\ref{sigm/ua})   we have
\beqn\label{pr/P1}
v_0 u_\as(R) =a(-R) P[\sigma_R<T]\{1+o(1)\} \geq - P(\La_R)/A(R)\{1+o(1)\}.
\eeqn
Here for the inequality we have employed the second condition of (\ref{e5.1}) in addition to  (\ref{L4.8}), (\ref{Dual}) and  $a(-x)-a(x)\sim -1/A(x)$. 
Thus the required lower bound of $u_\as$ is obtained. Since $a(-R)A(R)$ is bounded because of (b), we have $P_0[\sigma_R <T] \geq cP_0(\La_R)$ for some  $c>0$.  This shows the last formula of the proposition, for 
$P[\sigma_R <T] \sim P[\sigma_{R+1} <T] \leq P(\La_R)$. 
\qed

The following  lemma is used in the next subsection.
\begin{lemma}\label{lem4.4(2)} Suppose that  (\ref{e5.1}) holds and  $H_+$ is of dominated variation.
\v2
{\rm  (i)}\;\;  For any  $1/2\leq \de<1$, 
\beqn\label{eqL80}
g_\vom(w, y) 
\left\{\begin{array} {ll}
 =V_\ds(w)\times  o\big(U_\as(R)\big/R\big)  \quad &0\leq w<\de y,\\[1mm]
  \leq a(-y) \leq C/|A(y)| \quad &\de y\leq  w\leq y/\de,\\[1mm]
  \sim  U_\as(y)/\hat\ell^*(w) \quad &w>y/\de >0.
  \end{array}
  \right.
  \eeqn

{\rm  (ii)}\; \; ${\displaystyle \sup_{y\geq 0}\; \sum_{w=0}^\infty g_\vom(w,y)} [1-F(w)]<\infty$. 
\end{lemma}
\pf  Let  $EX=0$.  
By  (\ref{G-ft}), the Spitzer's representation of $g_\vom(w,y)$, one can easily deduce  (i),  
 with the help of  (\ref{u/U}) and (\ref{A/VU}). For convenience of later citations we note that (i) 
 entails  that
for some constant $C$
\beqn\label{eqL8}
g_\vom(w, y) 
\left\{\begin{array} {ll}
  \leq C V_\ds(w)U_\as(y)\big/y \quad & 0\leq w< y/\de,\\[1mm]
  \leq  C U_\as (y)/\hat\ell^*(w)  \quad &w\geq y/\de >0.
 \end{array} \right.
\eeqn
  
   By (\ref{eqL8}) $\sum_{w=0}^\infty g_\vom(w,y)[1-F(w)]$ is bounded above by a constant multiple of 
\beqn\label{sumVF}
 \frac{U_\as(y)}{y}\sum_{w=0}^{2y} V_\ds(w)H_+(w)+  U_\as (y) \sum_{w=2y}^\infty \frac1{\hat\ell^*(w)}H_+(w)
 \eeqn
for $y\geq x_0$. The first term above approaches zero as $y\to\infty$, for by {\bf L(3.3)}  $V_\ds(w)H_+(w) =o(1/U_\as(w))$,
while
 the second sum equals  $\hat \ell_\sharp(2y) \sim 1/U_\as (y)$.
 Thus  (ii) follows.  \qed

 \v2\v2
{\bf 4.2.} {\sc Asymptotic forms  of $u_\as $ and $P(\La_R)$.}
\v2

 Throughout  this subsection we suppose that  (\ref{e5.1})  hold. The results given in the preceding  subsection  except for Lemma \ref{lem4.4(2)} accordingly are applicable; in particular we  have the bound on $g_\vom(x,y)$ in the last lemma as well as the bounds
 \beqn\label{u/La}
 u_\as(y) =o\big(U_\as(y)\big/ y\big) \quad \mbox{and}\quad P(\La_y) \asymp P[\sigma_y<T] \sim v_0u_\as(y)/a(-y).
 \eeqn


   Recalling $B(R) = \mathbb{Z} \setminus [0,R]$, one sees
$$ P_x(\La_{R}) = \sum_{w=0}^{R} g_{B(R)}(x, R-w)H_+(w).$$
For each $r=1,2,\ldots$,  $g_{B(r)}(x,r-y), \, x,y \in B(r)$ is symmetric:
 \beqn\label{sym}
 g_{B(r)}(x,r-y) = g_{B(r)}(y, r-x), 
 \eeqn
 for the both sides equal $g_{- B(r)}(-r+y,-x)$.
Note that under (NRS), by duality,  $\hat Z$ is r.s. and 
 \beqn\label{eq4.4}
 v_\ds(x)\sim \frac1{\hat \ell^*(x)}\quad \mbox{and}\quad  U_\as   (x) \sim \frac1{\hat \ell_\sharp(x)} \quad \mbox{where}\quad \hat \ell_\sharp(t)= \int_t^\infty \frac{H_+(t)}{\hat\ell^*(t)}dt.  
 \eeqn
By Theorem \ref{thm1} and Remark \ref{rem3.1}(b) it follows that under  (\ref{e5.1}), 
 \beqn\label{A/VU}
 a(-x)-a(x)  \sim  \frac{-1}{A(x)} \sim \frac1{\hat \ell^*(x)\hat \ell_\sharp(x)} \sim  \frac{V_\ds(x)U_\as   (x)}{x}.\qquad
 \eeqn
  The next lemma is crucial for the proof of Theorem \ref{thm2}.   Recall that $N(R) = \sigma_{(R,\infty)} -1$ and 
that the condition  (\ref{ctn}) reads
$$ \limsup H_+(x/\lambda)/H_+(x) \to 1 
 \;\;\mbox{as} \;\lambda \downarrow 1.$$

\begin{lemma}\label{lem4.50}  Suppose that (\ref{e5.1})  holds and  $H_+$ varies dominatedly. Then for some constant $C$,
  \beqn\label{eq5.51}
 P[S_{N(R)} \geq {\textstyle \frac18}R, \, \La_R]\, \leq CU_\as(R)H_+(R) + o\big(P(\La_{R/2})\big).
 \eeqn 
 and if one further supposes the continuity condition (\ref{ctn}), 
 \beqn\label{eq5.52}
P[S_{N(R)} \geq {\textstyle \frac18}R,\, \La_R]\, \leq o\big(U_\as(R)H_+(R)\big) + o\big(P(\La_{R/2})\big).
 \eeqn 
\end{lemma}
\pf \, We use the representation
\beqn\label{start}
 P[S_{N(R)} \geq {\textstyle \frac18}R,  \La_R] = \sum_{0\leq w\leq R/8} g_{B(R)}(0,R-w)H_+(w).
 \eeqn

Splitting the r.w. paths by the landing points, $y$ say, when  $S$ started at the origin exits $[0,\frac12 R]$, we obtain for $0\leq w<R/8$ 
$$g_{B(R)}(0, R-w) 
=\sum_{R/2<y \leq R}P[S_{\sigma_{B(R/2)}}=y]  g_{B(R)}(y, R-w).$$
Hence
$$g_{B(R)}(0, R-w) 
= \sum_{z=0}^{R/2}\sum_{R/2<y \leq R}g_{B(\frac12 R)}(0,z)p(y-z) g_{B(R)}(y, R-w),$$
where $ p(x) = P[X=x]$.
 Taking  $0<\ep <1/8$ arbitrarily   we decompose the double sum into the following three parts:
$$I_\ep(w) = \sum_{z=0}^{R/2} \sum_{(1-\ep)R <y\leq R},  \quad I\!I_\ep(w) = \sum_{z=0}^{R/4} \sum_{\frac12 R<y\leq (1-\ep) R} \quad\mbox{and}\quad   I\!I\!I_\ep (w)= \sum_{\frac14 R < z \leq \frac12 R} \; \sum_{\frac12 R<y\leq \frac (1-\ep)R}. $$

First we evaluate $I\!I_\ep (w)$. 
Employing the symmetry (\ref{sym})  we see
\beqn\label{trv}
g_{B(R)}(y, R-w) = g_{B(R)}(w, R-y) \leq g_{\vom}(w, R-y).
\eeqn
Since  $g_\vom(w,R-y) \leq CV_\ds(w)U_\as(R)\big/(\ep R)$  in its range of $y$ by virtue of (\ref{eqL8}), we infer that 
\begin{eqnarray}
\frac{I\!I_\ep (w)}{U_\as (R)H_+(R)} &\leq&\frac{CV_\ds(w)\sum_{z=0}^{R/4} g_{B(\frac12 R)}(0,z)H_+({\textstyle \frac14} R)}{\ep RH_+(R)} 
 \nonumber\\
&\leq &\ep^{-1} C V_\ds(w)U_\as(R)\big/R.
\label{bd_II}
\end{eqnarray}
By {\bf L(3.3)}  $V_\ds(w)H_+(w)=o\big(1/U_\as(w)\big)$, so that 
\beqn\label{e4.26}
\int_0^R V_\ds(w)\mu_+(w)dw =o\big(R/U_\as(R)\big).
\eeqn
 Hence
\beqn\label{Ub_II}
\frac{1}{U_\as (R)H_+(R)}\int_0^R I\!I_\ep (w)H_+(w)dw \,\longrightarrow \, 0.
\eeqn

By the first case of (\ref{eqL80}), we have as above 
$$  I\!I\!I_\ep (w)= \bigg[V_\ds(w) \times o\bigg(\frac{U_\as(R)}{R}\bigg)\bigg]\sum_{\frac14 R<z\leq \frac12 R} \sum_{\frac12 R< y \leq\frac (1-\ep)R}g_{B(\frac12 R)}(0,z)p(y-z).
 $$ 
Making the change of variables  $z'= \frac12 R-z$ and $y'= y-\frac12 R$ shows that the double sum is less than
 $$\sum_{0\leq z'<\frac12 R}g_{B(\frac12 R)}(0, {\textstyle \frac12}R -z')H_+(z') = P(\La_{R/2}). 
$$
Hence applying  (\ref{e4.26})  we find 
\beqn\label{Ub_III}
\frac1{U_\as (R)H_+(R)}\int_{0}^{R}   I\!I\!I_\ep (w)H_+(w) dw = o\big(P(\La_{R/2})\big). 
\eeqn
  
Because of (\ref{Ub_II}) and (\ref{Ub_III}) as well as  Lemma \ref{lem4.2}(i),  for verification of (\ref{eq5.51})  it  suffices to show that  for an $\ep>0$ fixed,
\beqn\label{4.31}
\frac1{U_\as (R)H_+(R)}\sum_{w=0}^{R} I_ {\ep} (w)H_+(w) \leq C. 
\eeqn
Since $g_{B(R)}(0,z) < v_0u_\as (z)$, the above sum is less than 
$$ \sum_{w=0}^{R}   H_+(w) \sum_{z=0}^{R/2} \;  \sum_{y=(1-\ep) R}^R  g_{B(R)}(y,R-w) p(y-z)  u_\as (z).
 $$
Since $ \sum_{w=0}^{R}   H_+(w)g_{B(R)}(y,R-w) = P_{y}(\La_R)\leq 1$,   this triple sum is less than 
 \beqn\label{e5.31}
 \sum_{z=0}^{R/2} \; \sum_{y=(1-\ep) R}^R p(y-z)  u_\as (z).
\eeqn
It is clear that this double sum is at most a constant multiple of  $U_\as(R)H_+(R)$. Thus (\ref{4.31}) follows.
Since the  inner sum $\sum_{y=(1-\ep) R}^R p(y-z)$ equals
$ \mu_+(R-\ep R-z)- \mu_+(R-z-1)$, the continuity condition (\ref{ctn}) assures that    uniformly  for $z\leq R/2$,   as $\ep\to 0$  the above double sum is of the smaller order of magnitude than $U_\as (R)H_+(R)$ so that 
\beqn\label{4.32}
\frac1{U_\as (R)H_+(R)}\sum_{w=0}^{R} I_ {\ep} (w)H_+(w) \to 0 \quad \mbox{as\; $R\to\infty$\; and\; $\ep\downarrow 0$\; in this order}. 
\eeqn
showing  (\ref{eq5.52}).
 \qed
 
\begin{lemma}\label{lem4.5}  Suppose that (\ref{e5.1}) holds and $H_+$ varies dominatedly. 
  Then 
  \beqn\label{asy-u} 
P(\La_R) \asymp U_\as (R)H_+(R) \quad \mbox{and} \quad  u_\as (x) \asymp U_\as (x)H_+(x)a(-x),
\eeqn
and if one further supposes the continuity condition (\ref{ctn}), 
 \beqn\label{eq5.5}
 P[S_{N(R)} \geq {\textstyle \frac18}R\, |\, \La_R]\, \longrightarrow\, 0.
 \eeqn 
\end{lemma}
\pf   Lemma \ref{lem4.2}(ii) and \ref{lem4.50} together show
$$P(\La_R) \leq CU_\as (R)H_+(R) + o\big( P(\La_{R/2})\big).$$
On putting  $\lambda(R) =P(\La_R)/[CU_\as (R)H_+(R)]$,  the dominated variation of $H_+$ allows us to rewrite this inequality as
$$\lambda(R) \leq 1 + o\big(\lambda(R/2)\big).$$
One can choose $R_0>0$ so that  $\lambda(R) \leq 1 + \frac12\lambda(R/2)$ if $R\geq R_0$, which yields
$$\lambda(R) \leq 1 + 2^{-1}\big[1+ 2^{-1}\lambda(R/4)\big] \leq \ldots \leq 2 + 2^{-n}\lambda(2^{-n} R)$$
as far as  $2^{-n}R\geq R_0$. Taking $n=n(R)$ so that $2^{-n}R\geq R_0>2^{-n-1}R$, one obtains that  $P(\La_R) = O\big(U_\as (R)\mu_+(R)\big)$, which shows (\ref{asy-u}) (recall the second  of (\ref{u/La})). This entails $P(\La_{R/2}) \asymp P(\La_{R}) $, hence  (\ref{eq5.5}) follows from the second half of the preceding lemma. \qed

\begin{rem}\label{rem4.1} The continuity condition (\ref{ctn}), used at the end of  the proof of Lemma \ref{lem4.50},   is necessary for $P_0(\La_R)/v_0 \sim U_\as (R) \mu_+(R)$  to hold
[note that the contribution to  the sum (\ref{e5.31}) from $z <\ep R$ signifies for any $\ep>0$]. Also $P(\La_y)/ u_\as (y)$ is not asymptotic to  $v_0A(y)$  if (\ref{ctn}) fails to hold.
\end{rem}
 
 As a consequence of Lemmas \ref{lem4.2} and \ref{lem4.5}  we obtain
\begin{lemma}\label{lem4.6}   Suppose that  (\ref{e5.1}) and (\ref{ctn}) hold.  Then
\v2
{\rm  (i)}\;  for each $\ep >0$, as $R\to\infty$
\beqn\label{eqL10}
P[ S_{N(R)} \geq \ep R\,|\, \La_R]\to 0,
\eeqn
and  $P_0[ \ep R\leq Z(R)\leq \ep^{-1} R\,|\, \La_R] \to 0$ as $\ep \downarrow 0$ uniformly in $R>1\,;$
\v2
{\rm  (ii)}\quad\;
$P(\La_R)/v_0 \sim U_\as (R)H_+(R)$ \quad and \quad
${\displaystyle \frac{P[\sigma_R<T]}{P(\La_R)}   \sim \frac{a(-R)- a(R)}{a(-R)}   } \,;$

\v2\v2
{\rm  (iii)} \quad \;
$u_\as (y) \sim {\displaystyle \frac{ U_\as (y)H_+(y)}{-A(y)}.}$
\end{lemma}
\v2
\pf
The first convergence of (i) follows from Lemmas \ref{lem4.2}(ii) and  \ref{lem4.5}, and  the second one of  (i)  from it---by virtue of (\ref{ctn}). By  (\ref{eqL10}) we have
$$P(\La_R) \sim P[ S_{N(R)} \leq \ep R,  \La_R] 
\sim v_0 \int_0^{\ep R} u_\as (t) H_+(R-t)dt.$$
The integral of the last member is between $U_\as(\ep R)\mu_+(R-\ep R)$ and $U_\as(\ep R)\mu_+(R)$, and hence  
 may be written as $U_\as (R)\mu_+(R)\{1+o_\ep(1)\}$  as $R\to\infty$ and $\ep\downarrow 0$ owing to (\ref{ctn}) again. Since  $U_\as$ is s.v., this shows  the first relation of (ii). 
  By the second of  (i) 
$$P[\sigma_R<T] \sim P(\La_R)\bigg(o_\ep(1)+\sum_{y= \ep R}^{\infty} P[ Z(R) = y\,|\, \La_R]P_y[\sigma_R<T] \bigg),$$
where $o_\ep(1)\to 0$ as  $\ep \downarrow 0$. By (\ref{Dual}) the second probability under the summation sign is asymptotically equivalent to   $[a(-R)-a(R)]/a(-R)$, and we have the second relation of (ii). By (\ref{P_a}) 
we have $-1/A(y) \sim a(-y)-a(y)$, hence by the second equivalence of (ii)
$$P[\sigma_y<T] \sim   P(\La_y)/|A(y)a(-y)|,$$
while the probability on the LHS $\sim v_0u_\as (y)/a(-y)$.
 Combined with the first one of (ii) this yields (iii). \qed
 

\v2
 \v2
{\bf 4.3.}  {\sc  Proof of  Theorem \ref{thm2}.}

\v2
By virtue of Lemma \ref{lem4.6} Theorem \ref{thm2} follows if we can show the following
\begin{proposition}\label{prop4.1} Suppose that  (\ref{e5.1})  and  (\ref{ctn}) hold. 
Then  
 for each $\de<1$,  uniformly for  $0\leq x < \de y$,  as $y \to\infty$
\beqn\label{xLa00}
g_\vom(x,y)  \sim 
  \frac{V_\ds(x)U_\as (y)}{|A(y)|(x+1)} \sum_{k=0}^x H_+(y-k), \quad
P_x(\La_y) \sim g_\vom(x,y)A(y)
 \eeqn
and  
 $$ \forall \ep>0, \;\,  P_x[S_{N(y)} \geq x+\ep y\, |\, \La_y] \to 0.$$ 
\end{proposition}
\v2
\begin{rem}\label{rem4.3} The continuity condition (\ref{ctn}) is necessary for $P(\La_R) \sim  V_\ds(x)P(\La_R)/v_0$ to hold, as is seen from the identity $P(\La_y)/v_0 = \sum_{z=0}^R u_\as (z) H_+(R-z)$ where the contribution of the sum over $z <\ep y$ always signifies for any $\ep>0$. Also $P(\La_y)/ u_\as (y)$ is not asymptotic to  $v_0A(y)$  if (\ref{ctn}) fails to hold.
\end{rem}

\begin{lemma}\label{lem4.7}   Under the same assumption as in Proposition \ref{prop4.1} for each $\de<1$,  as $y\to\infty$
\beqn\label{gxy}
g_\vom(x,y)  \sim 
  \frac{U_\as (y)V_\ds(x)}{|A(y)|(x+1)} \sum_{k=0}^x H_+(y-k) 
  \quad \mbox{ uniformly for \;} 0\leq x < \de y.
 \eeqn
\end{lemma}
\pf  Substituting the asymptotic form  of $u_\as $ of Lemma \ref{lem4.6} into  (\ref{G-ft}) one has 
$$ g_\vom(x,y) \sim \frac{U_\as (y)}{-A(y)} \sum_{k=0}^x v_\ds(k) H_+(y-x+k) \quad \mbox{uniformly for \;} 0\leq x <\de y.$$ 
It accordingly is enough to see that the above sum may be replaced by $(x+1)^{-1}\sum_{k=0}^x\mu_+(y-k)$ as $x\to\infty$.  For each $x$ fixed, this is obvious  because of (\ref{ctn}). In case  $x\to\infty$ observe that the sum  restricted to $k<\ep x$ is less than $\ep x \mu_+(y-x)/\hat \ell^*(x)$  which is at most $C\ep$ of the remaining sum.
Then the result follows. 
\qed

\begin{lemma}\label{lem4.8}  Suppose  the same assumption as in  Proposition \ref{prop4.1} to hold. Then for each  $\de <1$, uniformly for $0\leq x<\de R$ 
 $$P_x[S_{N(R)} \geq {\textstyle \frac12} (1+\de) R\, |\, \La_R] \to 0.$$ 
\end{lemma}
\pf   
Put $R'= \lfloor \de R\rfloor$,  $R''=R-R'$ and $Q_{R}(x) = P_x[S_{N(R)}> R- \frac12 R'', \La_R]$.
Then
$$\frac{Q_{R}(x)}{P_x(\La_{R'})} = \sum_{k=R'}^{R} P_x\big[Z(R') =k\,\big|\, \La_{R'}\big] Q_{R}(k).$$ 
 Since for each $\ep>0$,  $Q_R(k) < P_k(\La_R) \to 0$ uniformly for $R'\leq k\leq R-\ep R''$ and $P_x(\La_{R'}) \asymp P_x[\sigma_R <T] \asymp P_x(\La_R)$ according to the second half of Lemma \ref{lem4.3}, we have only to show that uniformly for $0\leq x <R'$, as $\ep\downarrow 0$,
\beqn\label{L5.9}
\limsup_{R\to\infty} \sup_{0\leq x<R'} \frac1{P_x(\La_R)}\sum_{0\leq y <\ep R''} P_x\big[Z(R') =R-y\,\big|\, \La_{R'}\big] Q_{R}(R-y) \, \longrightarrow\, 0.
 \eeqn
As in {rem3.1(2)}the proof of Lemma \ref{lem4.5} we see that the sum above is dominated by 
\beqn\label{start1} 
 \sum_{0\leq y <\ep R''} \sum_{z=0}^{R''}g_\vom(x,z) p(R-y-z) \sum_{w=0}^{\infty}   g_\vom(w, y)H_+(w).
 \eeqn
which by Lemma \ref{lem4.4(2)}(ii) is at most a constant multiple of
$$  \sum_{z=0}^{R''}g_{\vom}(x,z) [F(R-z)- F(R-\ep R'' -z)].  $$
By {\bf L(2.1)}  $ \sum_{1\leq y <\ep R''} \sum_{z=0}^{R''}g_\vom(x,z) \leq V_\ds(x) U(R)$. Since  
$R-R'' \geq  \de R-1$, by (\ref{ctn}) the above sum is  of the smaller order of magnitude  than $V_\ds(x) U_\as (R) \times H_+(R) \asymp P_x(\La_R)$ as $R\to\infty$ and $\ep\downarrow 0$ in this order. 
Thus we have (\ref{L5.9}) as desired. \qed
 \v2
 
{\bf Proof of Proposition \ref{prop4.1} and Theorem \ref{thm2}.}  Since, by the asymptotic form of $u_\as$ in Lemma \ref{lem4.6}(iii),   $g_\vom(x,z) \leq CV_\ds(x)u(z)$ for $z\geq x+ \ep y$, we see $P_x[x+ \ep y< S_{N(y)} < (1-\ep) y\,|\, \La_y] \to 0$ uniformly for $0\leq x< \de y$ as in the proof of Lemma \ref{lem4.2}(ii), and the last formula of Proposition \ref{prop4.1} follows from Lemma \ref{lem4.8}.
 
As before we deduce from Lemma \ref{lem4.3}   and Lemma \ref{lem4.8} that uniformly for \, $0\leq x <\de y$
\beqn\label{cL10}
P_x[ \ep y < Z(y) < y/\ep\,|\,\La_y] \to 1 \quad \mbox{as \; $y\to\infty$ and  $\ep \downarrow 0$ in this order}.
\eeqn 
This shows that as $y\to\infty$
\beqn\label{e5.40}
 P_x(\La_y) \sim -A(y)a(-y)P_x[\sigma_{y} <T] \sim -A(y)g_\vom(x,y) \quad \mbox{uniformly for $0\leq x< \de y$.}
\eeqn
The first equivalence relation is the same as that giving the asymptotic form of $P_x[\sigma_{y} <T]$ in Theorem \ref{thm2},
and the second implies the asymptotic form of $P_x(\La_y)$ in (\ref{xLa00}). This finishes proof of Proposition \ref{prop4.1} (hence of  Theorem \ref{thm2}),  the other assertions being contained in Lemmas \ref{lem4.6} and \ref{lem4.7}. \qed

\section {Proof of Proposition \ref{prop1.1}  (transient case) and Theorem \ref{thm3}}

 Suppose that $F$ is transient and n.r.s.   According to \cite{Urenw},  we then have
    \beqn\label{Gpm}
 \begin{array}{ll}
 {\rm (a)} \;\; A(x)\to -\infty,\; G(-x) \;  \mbox{is s.v.   and} \;\; G(-x)- G(x) \sim -1/A(x). \\[2mm] 
 {\rm (b)}  \;\; 
   {\displaystyle G(-x) \sim \int_x^\infty \frac{F(-t)}{A^2(t)}dt,\;  G(x) \sim \int_x^\infty \frac{H_+(t)}{A^2(t)}dt + o(G(-x))} \qquad\qquad
 \end{array}
 \eeqn
 \v2\noindent
 as $x\to\infty$.
 (The last result in (a) is not stated in \cite{Urenw}, but actually proved in the proof of Theorem of \cite{Urenw} [see Eq(57), Eq(30)  and Section 3.3 of  \cite{Urenw}.)

 \begin{lemma}\label{lem5.1}  If $F$ is  transient and n.r.s., then $g_\vom(x,x)  \to G(0)$, and 
 for any $M>1$, as $y\to\infty$
 $$g_\vom(x,y) = G(y-x) -G(y) +o(G(-y)) \quad \; \mbox{uniformly for} \;\; 0\leq x< My.$$
 \end{lemma}
 \pf  Under (NRS) $\hat Z$ is r.s., so that $P_x[S_T<-\ep x] \to 0$ for each  $\ep>0$. Hence the result is immediate from the identity
$g_\vom(x,y) = G(y-x) - E_x[G(y-S_T)] $.
\qed 
\v2 
 
From Lemma \ref{lem5.1} one infers that as $y\to \infty$
\beqn\label{eL1}
\begin{array}{ll}
{\rm (i)}\;\; g_\vom(x,y) = \left\{
\begin{array}{lr}
G(-y) -G(y) + o(G(-y)) \quad \; \mbox{uniformly for} \;\; y/\de <x< My,\\
G(y-x)\{1+o(1)\}  + o(G(-y)) \;\;\;\;  \mbox{if\; $G(y)/G(y-x)\to 0$},\\
   o(G(-y)) \quad \; \mbox{uniformly for} \;\; 0\leq x< \de y; \,\mbox{and}
   \end{array}\right.\\[8mm]
{\rm (ii)}\;\; {\displaystyle \sum_{\frac12 y\leq w\leq \frac23 y} g_\vom(w,y) \leq O\big(yG(-y)\big).}
 \end{array}
 \eeqn 
 The estimate of $g_\vom(x,y)$  in (i) is not exact for $x \leq y/\de$ (unless  $|x-y|/y$ is sufficiently small), but  for $x/y\asymp 1$ (ii) provides a bound enough for our present purpose.  

\v2
{\bf Proof of Proposition \ref{prop1.1} (transient case}). Under 
 (\ref{e5.3}),  by Lemma \ref{lem4.3}  we have  $P[Z(R) >MR\,|\,\La_R] \to 0$ as $R\wedge M \to \infty$,  while  under (\ref{e1.7}),  by the first case of (\ref{eL1}),  
 $$P_z[\sigma_R < T]\geq  [G(-R)-G(R)]/G(0)\{1+o(1)\}\quad\mbox{for\;\;  $R<z <MR$, $M>1$.}$$
  It therefore follows that    
$P[\sigma_R<T\,|\, \La_R] \geq  [G(-R)-G(R)]/G(0)\{1+o(1)\}.$
  Hence 
  \beqn\label{pr/P1t}
v_0 u_\as(R) =G(0) P[\sigma_R<T]\{1+o(1)\} \geq - P(\La_R)/A(R)\{1+o(1)\}.
\eeqn
The lower bound (\ref{LaP1}) of  $P(\La_R)$ also is valid owing to Lemma \ref{lem4.2},  which together with (\ref{pr/P1t}) yields that of $u_\as(x)$ in Proposition \ref{prop1.1}.   \qed
\v2

For the proof of Theorem \ref{thm3} we need to obtain the bound 
\beqn\label{Ub_P}
P[\sigma_R<T\,|\, \La_R] \leq  CG(-R).
\eeqn
Using  $g_\vom(R+y,R) \leq G(-y)$ we see that
 $$P[\sigma_R<T\,|\, \La_R] =\sum_{y=1}^\infty P[Z(R) =y\,|\, \La_R] \frac{g_\vom(R+y,R)}{G(0)}
 \leq  \frac{E[G(-Z(R))\,|\, \La_R]}{G(0)}.$$
 Therefore, it follows that  for each $\ep>0$,
\beqn\label{lem5.11}  
P[\sigma_R<T\,|\, \La_R] < G(-R)\{1+o(1)\}  +  E\big[G(-Z(R)); Z(R)< \ep R\,\big|\, \La_R\big].
\eeqn
\v2
 \begin{lemma}\label{lem5.20}   If   (\ref{e5.3}) and (\ref{TRNS})  (given in (c) and (d) of (\ref{e5.12}) below) hold in addition to (NRS), then
 \v2
    {\rm (i)} \quad $P[Z=x] \sim H_+(x)/\hat\ell^*(x)$; and
 \v2
   {\rm (ii)} \quad $u_\as(y)\sim U_\as(y)H_+(y)/|A(y)|$.
   \v2
   \end{lemma}
   \pf \, Since $v_\ds(x)\sim \hat\ell^*(x)$,  (\ref{TRNS}) entails the bound $\sum_{y=0}^{\ep x} v_\ds(y)p(x+y) \leq C \ep v_\ds(\ep x)H_+(x)$ for  $\ep>0$ and $x$ large enough, and by using (\ref{e5.3}) it is easy to see (i).
   
  Because  (\ref{TRNS}) entails (\ref{ctn}), from (i)  one deduces  that 
 $$\limsup_{x\to\infty}\sup_{\de x< y <x}\big|P[Z=y]-P[Z=x]\big|/ P[Z=x] \to 0\quad\mbox{ as\;  $\de\downarrow 1$}.$$
   In Appendix we shall show that under this condition (ii) follows (see Lemma \ref{lem6.2.2}).
  \qed
\v2
 
\v2
In the rest of this section suppose that  the assumption of Theorem \ref{thm3} holds, namely  
\beqn\label{e5.12}
\left\{
\begin{array}{ll}
{\rm (a)}\;\; \mbox{$F$\;  is n.r.s.},\\[1mm] 
{\rm (b)}\;\; F \; \mbox{is transient} \quad \mbox{and}\quad  \limsup_{x\to \infty} {G(x)}/{G(-x)} <1,
\qquad\qquad\qquad\qquad\qquad
\\[1mm]  
{\rm (c)}\;\; 
\exists \lambda>1, \;\;  \limsup {H_+(\lambda t)}/{H_+(t)}  <1,\\[1mm]
{\rm (d)}\;\;  \exists C>0, \;\; p(x) \leq C\mu_+(x)/x.
\end{array}\right.
\eeqn  
\v2\noindent
\v2\noindent
 Note that (\ref{TRNS}) entails that  $H_+$ varies dominatedly and satisfies the continuity condition  (\ref{ctn}).

 \v2
 \begin{lemma}\label{lem5.2}  If (\ref{e5.12}) holds, then for some constant $C$,
 $$E[G(-Z(R))\,|\, \La_R]  \leq CG(-R).$$
   \end{lemma}
 
 \v2\noindent
  \pf  We may estimate $E[G(-Z(R)), Z(R\leq R\,|\, \La_R]$  which  is represented as
 $$\frac1{P(\La_R)} \sum_{w=0}^{R}\sum_{z=0}^{R} g_{B(R)}(0,R-w)p(w+z)G(-z).$$ 
By virtue of  (\ref{e5.12}d) (and the dominated variation of $\mu_+$) 
the outer sum restricted to  $w\geq R/5$ is at most a constant multiple of 
\beqn\label{G/G}
\; \frac1{P(\La_R)}  \sum_{w=R/5}^{R}g_{B(R)}(0,R-w)\frac{\mu_+(R)}{R}\sum_{z=0}^R G(-z) \leq G(-R)\{1+o(1)\}.
\eeqn
We show that the other sum is $o\big(G(-R)\big)$.
  To this end we proceed as in the proof pf Lemma \ref{lem4.5}. Let $I_\ep(w)$, $I\!I_\ep(w)$ and $I\!I\!I(w)$ be as therein. For the present purpose  we take $\ep=1/4$ and drop the subscript $\ep$ from   $I_\ep(w)$ and  $I\!I_\ep(w)$. 
  What is to be shown may then  be paraphrased as
  \beqn\label{eL5.2}
  \frac1{P(\La_R)}\sum_{w=0}^{R/8}\sum_{z=0}^{R} \big[I(w) +I\!I(w)+I\!I\!I(w)\big] p(w+z)G(-z) =  o\big(G(-R)\big).
  \eeqn
  First of all  we note that combining (\ref{e5.3}) and  (\ref{TRNS}) leads to 
  \begin{eqnarray}\label{HG}
  \sum_{z=0}^\infty p(w+z)G(-z) &\leq& C\frac{H_+(w)}{w}\sum_{z=0}^w G(-z) + H_+(w)G(-w)\nonumber\\
   &\leq&   C' H_+(w)G(-w) \quad (w\geq 1).
   \end{eqnarray}
   
  Recall 
   $$g_{B(R)}(0, R-w) 
=\sum_{R/2<y \leq R}g_{B(R/2)}(0,z)p(y-z) g_{B(R)}(y, R-w)$$
and that $I(w)$, $I\!I(w)$ and $I\!I\!I(w)$ are defined as the contributions to this sum from the ranges, 
   $$0\leq z \leq {\textstyle \frac12}R, \;{\textstyle \frac34} R < y\leq R; \quad 0\leq z \leq {\textstyle \frac14}R, \;{\textstyle \frac12} R < y\leq {\textstyle \frac34} R\quad \mbox{and} \quad {\textstyle \frac14}R < z\leq {\textstyle \frac12}R,\; {\textstyle \frac12}R<y \leq {\textstyle \frac34}R,$$
   respectively.  
 By  (\ref{TRNS})  it follows that 
   $$I(w) \leq \frac{CU_\as(R)H_+({\textstyle \frac14}R)}{R}\sum_{y=\frac34 R}^R g_\vom(w, R-y) \leq \frac{CU^2_\as(R)H_+(R)}{R}V_\ds(w),$$
where for the latter inequality  {\bf L(2.1)} is employed.
 Then, using  $P(\La_x) \geq U_\as(x)H_+(x)\{v_0+o(1)\}$,  (\ref{HG}) and $V_\ds(x)H_+(x)=o\big(1/U_\as(x)\big)$ in turn,  we obtain 
 $$\frac1{P(\La_R)}\sum_{w=0}^{R/8}\sum_{z=0}^R I(w) p(w+z)G(-z) \leq\frac{C'' U_\as(R)}{R}\sum_{w=0}^{R/8} V_\ds(w)H_+(w)G(-w) = o\big(G(-R)\big).$$
One can  obtain the corresponding bound for $I\!I(w)$ in the same way. 

 As for  $I\!I\!I(w)$, we apply Lemma \ref{lem5.20}(ii) to see that  
 \beq
  I\!I\!I(w) &\leq& \sum_{z=R/4}^{R/2} \sum_{y= R/2}^{3R/4}  u_\as(z)p(y-z)g_\om(w,R-y)\\
  &\leq& \bigg[C\frac{U_\as(R)H_+(R)}{A(R)}\bigg]^2 V_\ds(w)\sum_{z=R/4}^{R/2} \sum_{y= R/2}^{3R/4} p(y-z)\\
  &\leq & C'\bigg[\frac{U_\as(R)H_+(R)}{A(R)}\bigg]^2 V_\ds(w) A_+(R).
  \eeq
  As above we have $\sum_{w}^{R} U_\as(w) H_+(w)G(-w) <\!< RG(-R)/U_\as(R)$, so that
  $$\frac1{U_\as(R)H_+(R)}\sum_{w=0}^R  I\!I\!I(w)  \mu_+(w) G(-w) <\!< \frac{R H_+(R)G(-R)}{A(R)}= o\Big(G(-R)\Big).$$
Thus (\ref{eL5.2}) is verified as required. 
 \qed

\v2
\begin{lemma}\label{lem5.3}  Under the same assumption as Lemma \ref{lem5.2},
\v2
{\rm (i)}\;\, $E[G(-Z(R)), Z(R)<\ep R\,|\,\La_R]/G(-R) \to 0$ \;\, as\;$R\to\infty$\;and\;$\ep \downarrow 0$\;in this order;\, and
\v2
{\rm (ii)}\;\; $P[S_{N(R)} > \ep R\,|\,\La_R]\to 0$ \quad for each  $\ep>0$.
\end{lemma}
\v2\noindent
\pf  In the  proof of Lemma \ref{lem5.2}  we have seen that 
\beqn\label{eL5.3}
\frac1{P_0(\La_R)} \sum_{w=0}^{R/5}\sum_{z=0}^{R}  
g_\vom(0, R-w) \frac{\mu_+(w+z)}{w+z}G(-z) = o\big(G(-R)\big).
\eeqn
 On the other hand,   the trivial bound ${\mu_+(w+z)}/(w+z)\leq \mu_+(w)/w$  yields
$$ \frac1{P_0(\La_R)} \sum_{w=R/5}^R\sum_{z=0}^{\ep R}  g_{B(R)}(0, R-w) \frac{\mu_+(w+z)}{w+z}G(-z) \leq C \ep G(-R).$$
By the same reasoning as was advanced for the bound (\ref{G/G}) these together show (i).  On noting that  $\sum_{\frac15 R<w <R-\ep R} g_{B(R)}(0,R-w)\mu_+(w) =o(U(R)\mu_+(R))$, the first convergence of (ii) follows from (\ref{eL5.3}) since $G(-x)$ is asymptotically decreasing. \qed
\v2
By Lemmas \ref{lem5.3}(i) and \ref{lem4.3} $P[\ep R <Z(R)< \ep^{-1} R\,|\, |\La_R] \to 1$  as $R\to\infty$ and $\ep\downarrow 0$ in turn, while  for each $\de<1$,  
$$P_z[\sigma_R <T\,|\, \La_R] \sim [G(-R)-G(R)]/G(0)$$
 uniformly for  $0\leq z \leq \de R$.  With these relations together with Lemmas \ref{lem5.3}(i) again  the same arguments  as made  in case $EX=0$  lead to
\beqn\label{u/UH/A}
 P(\La_y) \sim U_\as(y)H_+(y)\quad \mbox{and}\quad u_\as(y)\sim P(\La_y)/|A(y)|.
 \eeqn

\v2
\begin{lemma}\label{lem5.5}  Under the same assumption as in Lemma \ref{lem5.2}, for each  $\de<1$,
\beqn\label{Ub_P_La}
P_x(\La_R) \leq C'U_\as(x) U_\as(R)H_+(R)\quad (0\leq x<\de R).
\eeqn
\end{lemma}
\pf\, 
Let $0\leq x<\de R$ for $\de<1$ and $\ep =\frac12 (1-\de)$.  By {\bf L(2.1)}
$$ \sum_{y=0}^{(1-\ep) R}
g_\vom(x, y)H_+(R-y) \leq U_\as(x) U_\as(R)H_+(\ep R).
$$
For $x\geq R/8$, using  the estimate of $u_\as$ in Lemma \ref{lem5.20} (ii) and {\bf L(3.3)} one infers 
\beq
 \sum_{w=0}^{\ep R}
g_{B(R)}(x, R-w)H_+(w)\leq   \sum_{w=0}^{\ep R}
g_{\vom}(w, R-x)H_+(w) &=& o\big(RH_+(R)/|A(R)|\big) \\
&=&o\big(U_\as(x) U_\as(R)H_+(R)\big).
\eeq
For $x<R/8$ one proceeds as in the proof of Lemma \ref{lem4.50} based on 
the representation
$$P_x(\La_R) = \sum_{w=0}^\infty \sum_{z=0}^{R/2}\sum_{y=R/2}^R  g_{B(R/2)}(x,z) p(y-z) g_{B(R)}(w,R-y) H_+(w).$$
The sum restricted to $z\geq \frac14 R$, $y\leq frac34 R$ is at most a constant multiple of 
$$U_\as(x) \bigg[\frac{U_\as(R)H_+(R)}{A(R)}\bigg]^2 RH_+(R) \frac{R}{U_\as(R)} =o\big(V_\ds(x) U_\as(R)H_+(R)\big).$$
The rest of the sum is readily evaluated to be $o\big(V_\ds(x) U_\as(R)H_+(R)\big)$. Thus (\ref{Ub_P_La}) is verified. 
\qed

\begin{lemma}\label{lem5.6}  Under the same assumption as in Lemma \ref{lem5.2}, for each  $\de<1$,
as\;$R\to\infty$\; and\;$\ep \downarrow 0$\; in this order
\v2
{\rm (i)} \quad $E_x[G(-Z(R)), Z(R)<\ep R\,|\,\La_R]/G(-R) \to 0$;\; and
\v2
{\rm (ii)} \quad $P_x[S_{N(R)} > {\textstyle \frac12} (1+\de)  R\,|\,\La_R]\to 0$\quad and \quad  $P_x[Z(R)<\ep R\,|\,\La_R]\to 0,$
 \v2\noindent
in both {\rm (i)} and {\rm (ii)}  the convergence being uniform for $0\leq x<\de R$. 
\end{lemma}
\pf  By the asymptotic form of  $u_\as$ obtained in  (\ref{u/UH/A}) we have the same asymptotic form  of
 $g_\vom(x,y)$ as given in Lemma \ref{lem4.7}, and this together with the  bound of $P_x(\La_R)$ in the preceding lemma allows us to follow
 the proof  of Lemma \ref{lem4.8} to show both (i) and (ii). \qed 
\v2
{\bf Proof of Theorem \ref{thm3}.} \, By virtue of  Lemma \ref{lem5.6}(i) we obtain
\beqn\label{5.40}
 P_x(\La_R)\sim \frac{G(0)P_x[\sigma_R <T]}{G(-R)-G(R)}   \sim -A(R)g_\vom(x,R) \quad \mbox{uniformly for $0\leq x< \de R$.}
\eeqn
With this as well as (\ref{u/UH/A}) 
we can follow the proof of Theorem \ref{thm2} to show the rest of the results of Theorem \ref{thm3}. 
\qed

\v2\v2
\section{Estimation of $P\big[ S_{N(R)}= y \,\big|\, \La_R\big]$}

In this section, we suppose 
\v2
$(*)$ \quad {\it either the  assumption of  Theorem \ref{thm2} or that of Theorem \ref{thm3} holds}
\v2\noindent
and  compute the conditional probability of $S$ exiting 
$B(R)$ through $y\in B(R)$, given $\La_R$ and $S_0=x\in B(R)$.  Denote it by $q_R(x,y)$:
$$q_R(x,y) = P_x\big[ S_{N(R)}= y \,\big|\, \La_R\big].$$
Let  $1/2< \de <1$. In Sections 4 and 5 we have shown that uniformly for $0\leq x<\de R$,
\beqn\label{4-5.1}
P_x[ \ep R < Z(R) < R/\ep\,|\,\La_R] \to 1 \quad \mbox{as \; $R\to\infty$ and  $\ep \downarrow 0$ in this order};
\eeqn
uniformly for $ 0\leq x < \de y$, as $y\to\infty$ 
\beqn\label{4-5.2} 
 P(\La_y)\sim -u_\as(y)A(y)\sim U_\as(y)H_+(y), \quad   P_x(\La_y)  \sim -A(y)g_\vom(x,y), 
\eeqn
\beqn\label{4-5.3}
g_\vom(x,y)  \sim 
  \frac{U_\as (y)V_\ds(x)}{|A(y)|(x+1)} \sum_{k=0}^x H_+(y-k);
\eeqn
and for any $0<\ep<1/2$,
\beqn\label{4-5.4}
\sum_{y: |x-y|<\ep x} g_\vom(x,y) \leq C\frac{\ep x}{|A(x)|}.
\eeqn
Because of the slow variation of  $v_\ds(z)\sim 1/\hat\ell^*(y)$ we also have   
\beqn\label{g/U}
g_\vom(x,y) \sim U_\as(y)\big/\hat \ell^*(x) \quad  0\leq y <\de x. 
\eeqn

\begin{lemma}\label{lem6.1} \, Under $(*)$,  uniformly for $0\leq x<\de R$, $0\leq y< R$,  as $R\to\infty$ 
\beqn\label{g_B}
g_{B(R)}(x,y) = g_\vom(x,y)  -  \frac {U_\as (y)}{U_\as (R)}\, g_\vom(x,R)\{1+o(1)\};
\eeqn
in particular uniformly for\; $0\leq y \leq R$, as $R\to\infty$
\beqn\label{g_B0}
g_{B(R)}(0,y)= v_0u_\as (y)\bigg[1-\frac{u_\as (R)/U_\as (R)}{u_\as(y)/U_\as (y)}\{1+o(1)\}\bigg].
\eeqn  
 \end{lemma}
\pf 
By  (\ref{g/U}) and   (\ref{4-5.1}) we  deduce  that uniformly for $0\leq x<\de R$, $0\leq y< R$,
$$E_x\big[g_\vom(S_{\sigma(R,\infty)},y)\,\big|\,\La_R\big]\sim  U_\as (y)/\hat \ell^*(R),$$
and then,  using $-A(R)/\hat\ell^*(R) \sim 1/U_\as (R)$ and (\ref{4-5.2}),  that
$$E_x\big[g_\vom(Z(R),y), \,\La_R\big]\sim \frac{g_\vom(x,R)U_\as (y)}{U_\as (R)} \quad \mbox{uniformly for $0\leq y<R, \, 0\leq x <\de R$}. 
$$
Thus (\ref{g_B}) follows. By $g_\vom(0,y)=v_0u_\as (y)$ (\ref{g_B0}) is immediate from (\ref{g_B}).  \qed
\v2

By Lemma \ref{lem6.1} and (\ref{4-5.2})
\beqn\label{qR}
q_R(x,y)   =\frac{g_{B(R)}(x,y)}{P_x(\La_R)}H_+(R-y) = \bigg[\frac{g_\vom(x,y)}{g_\vom(x,R)}\{1+o(1)\} - \frac{U_\as (y)}{U_\as (R)}\bigg]\frac{H_+(R-y)}{A(R)}.
\eeqn
We begin with the case $x=0$.  
\begin{lemma}\label{lem6.2}  Under $(*)$,    as $R\to\infty$
\beqn\label{La0}
q_R(0,y) =
\left\{\begin{array}{lr}
{\displaystyle \frac{u_\as (y)}{U_\as (R)}\{1+o(1)\} }
\quad  \mbox{as} \;\; y/R \downarrow 0, \\[4mm]
o\big(H_+(R-y)\big/A(R)\big) \qquad  \mbox{as}\;\;\;   y/R \uparrow 1,
\end{array}\right.
\eeqn
and if\, $xH_+(x)$ is s.v. in addition, then   for each $\ep>0$, the first formula of (\ref{La0}) holds uniformly for \, $0 \leq y<(1-\ep)R$.
\end{lemma}
\pf
The formula (\ref{La0}) follows from (\ref{g_B0}). Indeed, as $y/R\uparrow 1$  the latter entails 
$g_{B(R)}(0,y) = o(u_\as (R))$ while $u_\as (R) \sim P(\La_R)/A(R)$ in view of (\ref{4-5.2}). Hence,  the second case of  (\ref{La0}) is immediate from the identity
\beqn\label{ID}
q_R(0,y) =\frac{ g_{B(R)}(0,y) H_+(R-y)}{P(\La_R)}.
\eeqn
 In case $y/R\to 0$, it follows from (\ref{g_B}) that the RHS of the above identity is expressed as  $u_\as (y)/U_\as (R)\{1+o(1)\}$.
If  $xH_+(x)$ is s.v., then for  $\ep R <y<(1-\ep) R$,
 (\ref{g_B0}) entails  
 $$g_{B(R)}(0,y) =v_0u_\as (y)\big[1-y/R+o(1)\big],$$
  and hence
using  (\ref{4-5.2}) as well as (\ref{ID})  we see 
$$\qquad \qquad q_R(0,y)  \sim \frac{g_{B(R)}(0,y)}{v_0U_\as (R)}\cdot \frac1{1-y/R}\sim \frac{u_\as (y)}{U_\as (y)}. \qquad \qquad\qed  $$


\v2
Taking $0<\ep< 1/2$, we put  $\de =1-\ep$ and  
  let $x, y$ be such that $0\leq x\wedge y \leq x\vee y< \de R$ throughout the sequel. 
Substituting  (\ref{4-5.3}) and/or (\ref{g/U}) into (\ref{qR}) we obtain the following. 
  \v2
(a)\,  If $x\vee y <1/\ep$, then  
  $$q_R(x,y) \sim \frac{g_\vom(x,y)}{V_\ds(x)U_\as (R)}.$$  
  
(b)\,   Uniformly for  $y<\de x$,  as $x\to\infty$
  $$g_{B(R)}(x,y) \sim g_\vom(x,y) \sim V_\ds(x)U_\as(y)/x, \,\mbox{and}$$
  $$q_R(x,y) \sim \frac{U_\as (y)H_+(R-y)}{U_\as (R)\int_{R-x-1}^RH_+(t)dt} \asymp \frac{U_\as (y)}{U_\as (R)x}.$$  
  
(c)\,   Uniformly for $y>(x+1)/\de$, as  $R\to\infty$
 $$q_R(x,y) \sim 
 \left\{
 \begin{array}{ll}
 {\displaystyle \frac{H_+(R-y)}{A(R)}\bigg( \frac{\int_{y-x-1}^y H_+(t)dt}{\int_{R-x-1}^R H_+(t)dt}-1+o(1)\bigg) =o\bigg(\frac1{R}\bigg)  }\quad &\mbox{for} \;\; y>\ep R,\\[6mm]
 {\displaystyle \frac{u_\as (y)}{U_\as (R)} \cdot \frac{\int_{y-x-1}^y H_+(t)dt}{ (x+1)H_+(y)} 
  \asymp\frac{u_\as (y)}{U_\as (R)}  }\quad &\mbox{as} \;\; y/R \to 0.
  \end{array}\right.
  $$
  
 
\v2 \noindent
Using these estimates and (\ref{4-5.4}) we infer that  as $R\to\infty$
  \beqn\label{q_R}
  q_R(x,y)\left\{
   \begin{array}{ll}
   = o\big(1/R\big) \quad &\mbox{for} \;\; y> x/\de,\,  y>\ep R,\\[1mm]
   \sim  {\displaystyle \frac{u_\as (y)}{U_\as (R)} } \quad &\mbox{if} \;\; x/y\to 0,\,  y/R\to 0, \\[4mm]
   \geq  {\displaystyle \frac{U_\as (y)}{U_\as (R)x} }\{1+o_\de(1)\}  \quad &\mbox{if} \;\; \de x \leq y\leq x =o(R), \\[4mm]
 {\displaystyle  \sim \frac{U_\as (y)H_+(R-y)}{U_\as (R)\int_{R-x-1}^RH_+(t)dt} } \quad &\mbox{if} \;\; y<\de x, \, x\to\infty,
    \end{array}\right.
 \eeqn
\beqn\label{sum_q}
  \sum_{y:|x-y|< \ep' x}q_R(x,y) \leq  \frac{C\ep' x}{A(x)V_\ds(x)U_\as(R)} \sim 
 \ep' \frac{CU_\as(x)}{U_\as(R)}  \quad (0<\ep'< \ep/2),
 \eeqn
 Let $n_R$ be any function such that $U_\as  (n_R)\sim U_\as (R)$ and $n_R/R\to 0$. Then, on employing (\ref{sum_q}) and the first case of (c),
 \beqn\label{sum_q2}
   \sum_{y=x\vee n_R}^{R}q_R(x,y) \leq  
    \sum_{y=x\vee n_R}^{x\vee n_R + n_R}q_R(x,y) +  \frac{C \sum_{y=n_R}^{\de R} u_\as(y)}{U_\as(R)} +o(1) \,\longrightarrow\, 0,
 \eeqn
and we see that the mass of  the conditional distribution $P_x\big[ S_{N(R)} \in \cdot \,\big|\, \La_R\big]$ tends to concentrate on the set of $y$ such that
 \[
   \begin{array}{ll}
 \;\;  (1+x)/\ep <y< n_R \quad  &\mbox{if \; $U_\as (x)/U_\as (R) \to 0$, and}\\[2mm]
 \;\; o(x)< y < x,    &\mbox{if \;  $U_\as  (x)/U_\as (R) \to 1$},
  \end{array}
 \]
\v2\noindent
where  $o(x)$ is any function such that  $o(x)/x\to 0$ and $o(x)\to\infty$.  This suggests that the conditional walk moves to the right in the former case and to the left in the latter up to the epoch of exiting $B(R)$.  
 If $\ep <U_\as  (x)/U_\as (R) <\de$, the mass may be possibly distributed on both sides of $x$. 
 \v2
 {\bf Proof of Proposition \ref{prop1.3}.}  Let $\ep_0>0$ and $n_R$ be as above. Observe
 \[
 P_x\big[ Z(R) =z\, \big|\,\La_R \big] = \sum_{y=0}^R \frac{g_{B(R)}(x,y)p(z+R-y)}{P_x(\La_R)}= \sum_{y=0}^R q_{R}(x,y)\frac{p(z+R-y)}{H_+(R-y)}.
 \]
Then by  (\ref{sum_q2}) we see that uniformly  for $x > n_R$ and $z>\ep_0 R$, the last sum restricted to $y> x$ is negligible so that 
$$P_x\big[ Z(R) \leq z\, \big|\,\La_R \big]  = \sum_{y=0}^{x} \frac{g_\vom(x,y)}{P_x(\La_R)}P[R- y<X\leq R-y+z] + o(1). 
$$
We substitute the asymptotic form of $g_\vom(x,y)$ and $P_x(\La_R)$.  Noting (\ref{sum_q}) which allows us to replace  $g_\vom(x,y)$ by $U_\as(y)/\hat\ell^*(x)$, we infer that
the above sum is asymptotically equivalent to
$$\frac1{U_\as(R)\int_{R-x}^R H_+(t)dt}\sum_{y=0}^x U_\as(y)\big[H_+(R-y) - H_+(R+z-y)\big].$$
Now assuming $H_+(t)\sim L_+(t)/t$, we see that 
$$\sum_{y=0}^x U_\as(y) H_+(R-y)\sim  U_\as(x)\int_{R-x}^R H_+(t)dt \sim  -U_\as(R)L_+(R) \log\big[1-R^{-1}x\big],$$
and similarly $\sum_{y=0}^x U_\as(z) H_+(R+z-y) \sim -U_\as(R)L_+(R) \log[1-(R+z)^{-1}x].$ Since the ratio  $\log[1-(R+z)^{-1}x]$ to $\log[1-R^{-1}x]$ is  bounded away from 1 we obtain the formula of Proposition \ref{prop1.3} (with $z$ in place of $y$) in case $x> n_R$.  

For  $x\leq n_R$,  by  (\ref{q_R}) and (\ref{sum_q2}) again,  we see
 $$P_x\big[ Z(R) \leq z\, \big|\,\La_R \big]  = \sum_{y=0}^{2n_R} \frac{g_\vom(x,y)}{P_x(\La_R)}P[R- y<X\leq R-y+z] + o(1). 
$$
It is easy to see that $\sum_{y=0}^{2n_R} g_\vom(x,y) \sim V_\ds(x)U_\as(2n_R)$ (cf. \cite[Lemma 2.1]{Uexit}).  Since, by (\ref{ctn}), $H_+(R-y) - H_+(R+z-y)= H_+(R)- H_+(R+z) +o(H_+(R))$ for $0\leq y\leq x$,  we obtain 
$$P_x\big[ Z(R) \leq z\, \big|\,\La_R \big] = \frac1{H_+(R)}\big[H_+(R) - H_+(R+z)\big] +o(1).$$
Now the asserted formula of the proposition follows immediately. \qed

 \begin{rem}\label{rem6.1}  We can easily obtain the estimates of $P_x[S_{N(R)}=y \,|\,T< \sigma_{(R,\infty)}]$ corresponding to the results given above for the conditioning on $T<\sigma_{(R,\infty)}$ instead of $\La_R$ in  case  (NRS) with $\limsup a(x)/a(-x)<1$ or $\limsup G(x)/G(-x)<1$.    We consider it in the dual form.  Let $0<\ep<1/2 <\de<1$ and $0\leq x <\de R$.   By {\bf L(2.1)}  we have 
 $\sum_{y=0}^{\de R}  g_\vom(x,y)H_+(R-y) \leq V_\ds(x) U_\as (\de R) H_+((1-\de)R)$
and under (\ref{fun_res}) 
 $$P_x[S_{N(R)}<\de R\,|\,\La_R] \leq \frac{\sum_{y=0}^{\de R}  g_\vom(x,y)H_+(R-y)}{P_x(\La_R)}  \leq  V_\ds(R)U_\as (R) H_+((1-\de)R).$$
Owing to  {\bf L(3.3)} 
 it, therefore, follows that if  (C3) holds, then
 $$P_x[S_{N(R)}<\de R\,|\,\La_R] \, \longrightarrow\,  0,$$
 saying that the conditional distribution of $S_{N(R)}$ tends to concentrate in an interval contained in $[\de R, R]$ for any $\de <1$. (According to \cite[Theorem 1]{Uexit} the same is true under  the condition (AS) with $\alpha=2$ or $m_+(x)/m_-(x)\to 0$ where $m_\pm(x)= \int_0^x \eta_\pm(t)dt$.)

 Now suppose that (PRS) holds  and  that $\limsup a(-x)/a(x)<1$ or $\limsup G(-x)/G(x)<1$ according as  $EX=0$ or $E|X|=\infty$. Then  $A(x)\sim \ell^*(x)\ell_\sharp(x)$ and $g_\vom(x,y) = o(1/A(x))$ for  $x<y/\de$ (by the second formula of (\ref{g/V/a}) of Theorem \ref{thm1} and the dual of Lemma \ref{lem5.1}).  For $0\leq w<\ep R$ and $x>\de R$ one observes that
 $$E_w[g_\vom(S_{\sigma(R,\infty)},R-x);\La_R] \leq  P_w(\La_R)\sup_{y>R}g_\vom(y,R-x). $$
Since $P_w(\La_R) \sim V_\ds(w)/V_\ds(R)$, the RHS above is $o(V_\ds(w)/\ell^*(R))$, hence by the identity (\ref{idn_g}) one concludes that uniformly for  $0\leq w<\ep R <x\leq \de R$,
 \beqn\label{eR6.1}
 g_{B(R)}(x,R-w) \sim g_\vom(w,R-x)\sim V_\ds(w)/\ell^*(R).
 \eeqn 
so that 
 \begin{eqnarray}\label{eR6.11}
 P_x[S_{N(R)} =y \,|\,\La_R] &=& g_{B(R)}(x,y)\mu_+(R-y)/P_x(\La_R) \nonumber\\ 
 &=& \frac{1}{\ell^*(R)} V_\ds(w)\mu_+(w), \quad   w:=R-y.
 \end{eqnarray}
 Note that the sum of the last member over  $0\leq w < n_R$ tends to unity whenever $\ell^*(n_R)\sim \ell^*(R)$.  
 If one further supposes that $p(x) =O(\mu_+(x)/x)$ $(x\to\infty)$, then (\ref{eR6.1}) can be extended to $0\leq x <\de R$ where 
 $$g_{B(R)}(x, R-w) \sim P_x(\La_R)V_\ds(w)/\ell^*(R),$$
 and  (\ref{eR6.11}), accordingly,  holds  uniformly for $ 0\leq w< \ep R, 0\leq x <\de R$.
\end{rem}

\section{Proof of Proposition \ref{prop1.4}} 

We divide this section into two subsections.  In the first one we suppose  (AS) to hold and prove Proposition \ref{prop1.4} except for a part of  the last assertion (iv) of Proposition \ref{prop1.4}, which is dealt with in the second one under an  assumption  less restrictive than (AS). 
\v2\v2
{\bf 7.1.} {\sc Estimates of $g_\vom(x,y)$ under (AS). }
\v2
In this subsection we suppose (AS) to hold. Under (AS), there exist s.v. functions $\ell$ and $\hat \ell$ such that
\beqn\label{5.30}
U_\as   (x) \sim  x^{\al\rho}/\ell(x) \quad and \quad V_\ds(x)\sim  x^{\al\hat\rho}/\hat\ell(x)
\eeqn
as mentioned at the beginning of Section 5 of \cite{Uexit}, while it is shown in \cite{Uexit}(Lemma 5.1) 
$$U_\as(x)V_\ds(x)H(x) \,\longrightarrow\,   \big[q\al\rho B(\alpha\rho,\alpha\hat\rho)\big]^{-1}(\pi\al\hat\rho)^{-1}{\sin \pi\alpha\hat\rho} \quad\;\; (q>0),$$
where  $B(s,t) =\Ga(s+t)/\Ga(s)\Ga(t)$ with the understanding  $sB(s,t)=1$ if $s=0, t>0$.
By (\ref{U/Z})  we have  
$P[-\hat Z\geq x]/v_0 \sim (\pi\alpha\hat\rho)^{-1} (\sin \pi\alpha\hat\rho)/V_d(x)$ ($=o(1/V_\ds(x)$) if $\alpha\hat\rho=1$),   and these together with Lemma 4.6(ii) (in case $\alpha =q=\hat \rho =1$) yields
\beqn\label{hat_Z2}
P[-\hat Z\geq x]/v_0 \sim q\al\rho B(\alpha\rho,\alpha\hat \rho) x^{-\al\hat\rho\,}L(x)/\ell(x)
\eeqn
unless $q=\hat \rho=0$ (see Lemma \ref{lem6.1}). Recall that we may and do take $\ell= \ell^*$ or $\hat \ell_\sharp$ according as $\alpha \rho=1$ or $0$, and analogously for $\hat \ell$.

\begin{lemma} \label{lem*1} \, Suppose either  $\al=1$ with  $\rho\notin \{0,\frac12,1\}$
(entailing $p=q$)  or $1<\al<2$. Then 
\beqn\label{6.1}
{\rm  (a)}\;\; u_\as   (x) \sim \al\rho\, {x^{\alpha\rho-1}}/{\ell(x)}
\quad \mbox{and} \quad {\rm  (b)}\;\; v_\ds(x) \sim \al\hat\rho\, {x^{\alpha\hat\rho-1}}/\hat\ell(x). 
\eeqn
\end{lemma}
\pf 
 We prove only (b), (a) being  dealt with in the same way. 
First of all we recall that if $\al\hat\rho =1$, then $\hat Z$ is r.s.\,and the equivalence (b) follows (cf. 
\cite[Appendix B]{Upot}, \cite{Urenw}).
 It is  also noted that
in case  $1/2< \alpha\hat\rho <1$ the strong renewal theorem holds for $V_\ds$ without any extra assumption (cf. e.g., \cite{BGT}) so that (b) follows immediately from (\ref{5.30}). 


The proof for $\al\hat\rho\leq 1/2$ rests on the recent result by Caravenna and Doney \cite{CD}.
 According to Theorem 1.4 of \cite{CD}  it   suffices to show that if $\alpha\hat \rho\leq 1/2$,
\beqn\label{eq5.12}
\lim_{\ep\downarrow 0} \limsup_{x\to\infty} \, xP[-\hat Z \geq x]\sum_{z=1}^{\ep x} \frac{P[-\hat Z = x-z]}{z(P[-\hat Z \geq z])^2}=0.
\eeqn
Note that  (a)---as well as  (\ref{hat_Z2})---is applicable  since  
$\al \rho>1/2$ that entails  $q>0$. 
Writing  $p(\cdot)$ for  $P[X=\cdot]$,  we have the identities (equivalent to each other);   
\[ \frac{P[- \hat Z =   x]}{v_0} = \sum_{z=0}^\infty  u_\as   (z) p(-x-z),  \quad  \frac{P[- \hat Z \geq   x]}{v_0} = \sum_{z=0}^\infty  u_\as   (z) F(-x-z)
\]
 (see e.g. \cite[Eq(XII.3.6a)]{F}).
  We accordingly deduce that the sum in (\ref{eq5.12}) is dominated by a constant multiple of
$$J:=\sum_{z=1}^{\ep x} \sum_{y=1}^\infty \frac{ y^{\alpha\rho-1}}{\ell(y)} p(-x-y+z)
\frac{1}{z}
\bigg(\frac{z^{\alpha\hat\rho\,}\ell(z)}{L_-(z)}\bigg)^2,$$
where $L_-(xt) = F(-t)/t^\alpha$.
We may suppose  $y^{\al\rho-1}/\ell(y)$ to be decreasing. (If $\al\rho=1$, 
one may take $\int_0^x P[Z>t]dt/v_0$  for $\ell(x)$.)  Then we perform summation  by parts for the inner sum and, after replacing $F(-t)$ which thereby comes up by $qL(t)/t^\al$ with $L$ appropriately  chosen,
make summation by parts back as before to  obtain 
$$\sum_{y=1}^\infty \frac{ y^{\alpha\rho-1}}{\ell(y)} p(-x-y+z) \sim \al \sum_{y=1}^\infty \frac{qL(x+y-z)y^{\al\rho-1}}{\ell(y)(x+y-z)^{\al+1}} \leq C\frac{L(x)x^{\al\rho-\al-1}}{\ell(x)}.
$$
Hence
$$ J\leq C\frac{L_-(x)x^{\al\rho-\al-1}}{\ell(x)} \sum_{z=1}^{\ep x}  \frac{1}{z}
\bigg(\frac{z^{\alpha\hat\rho\,}\ell(z)}{qL(z)}\bigg)^2 \leq  C'\frac{\ell(x)}{L(x)}\ep^{2\al\hat\rho}x^{\al\hat\rho-1}.$$
Thus  by (\ref{hat_Z2}) $ xP[-\hat Z \geq x] J  \leq C'' \ep^{2\al\hat\rho}$, verifying (\ref{eq5.12}).
\qed

\v2
\begin{rem}\label{rem7.1} (a) The proof above depends on the fact that if 
$(\rho\vee\hat\rho)\al>1/2$, either  $Z$ or $\hat Z$ admits the strong renewal theorem. For this reason the case  $\al=2\rho=1$ is excluded from Lemma \ref{lem*1}, while the case $1/2 <\al<1$ may be included in it if $\rho$ satisfies the above condition. Any way,  taking the first formula of (\ref{T2-1}) in Theorem \ref{thm2} into account, we have the strong renewal theorem for $U_\as$ and $V_\ds$ at least in case $\alpha\geq 1$ under (AS) except in a few special cases that are (1) $\alpha=2\rho=1$; (2) $\alpha=\rho\vee \hat\rho=2p =1$; (3)  $\alpha=\rho\vee \hat\rho=1$ and $E|X|=\infty$,  of which  a partial result for the cases (2) and (3) is given by Lemma \ref{lem5.20} that entails  the following result: under (AS)
$$\frac{xu_\as(x)}{U_\as(x)}  \sim U_\as(x) V_\ds(x)\mu_+(x) \;\;\;\; \mbox{if}\;\; \alpha=\hat \rho=1\;\;  \mbox{and}\;\;  \limsup_{x\to\infty} \frac{p(x)|x|}{\mu_+(x)} <\infty;\;\mbox{and}$$
$$\frac{xv_\ds(x)}{V_\ds(x)}  \sim U_\as(x) V_\ds(x) \mu_-(x) \;\;\;\; \mbox{if}\;\; \alpha= \rho=1\;\;  \mbox{and}\;\;  \limsup_{x\to\infty} \frac{p(-x)|x|}{\mu_-(x)} <\infty.$$
\v2
(b) In the proof of Lemma \ref{lem*1} the property of the positive tail of $F$ is used only through those of the distributions of $Z$ and $\hat Z$. Since the regular variation of $u_\as$  and $\mu_-$ implies that of $P[-\hat Z>\cdot\,]$,  it accordingly follows---whether (AS) is true or not---that
\v2

{\it If $U_\as   (x)\sim x^{\beta}/\ell(x)$, $F(-x) \sim L_-(x)x^{-\al}$ with $L_-$ and $\ell$  s.v. and $0<\al-\beta<1$,   

 and $1/2<\beta \leq 1$, then (\ref{6.1}b) holds with $\al\hat\rho =\alpha - \beta$ and}
$$\hat \ell(x) =\frac{\Ga(\al-\beta+1)\Ga(\beta+1)}{\Ga(\al)\pi^{-1}\sin\al\hat\rho\pi}\cdot\frac{L_-(x)}{\ell(x)}.
$$

{\it If $Z$ is r.s.\,in particular, then from the condition $F(-x) \sim L_-(x)x^{-\al}$, $1<\al<2$, 
 
 it follows that $v_\ds(x) \sim \big[(\al-1)\pi^{-1}\sin (\al-1)\pi \big] x^{\al-2}\ell(x)/L_-(x)$. }
\v2
\end{rem}

Lemma \ref{lem*1} allows us to compute the precise asymptotic form of $g_{\vom}(x,y)$ for $\alpha\geq 1$ unless $\rho\in \{0,\frac12,1\}$, which case however is covered by Proposition \ref{prop3.1}.
 Note that   $g_{\vom}(0,y)=v_0u_\as(y)$ and $g_{\vom}(y,0) =v_\ds(y)$ for  $y\geq 0$.
\begin{lemma} \label{lem*2} \,
{\rm  (i)}  If $1<\al\leq 2$, 
\beqn\label{asymp_g}
g_{\vom}(x,y) \sim  \left\{ \begin{array} {ll} 
{\displaystyle \frac{\al\rho\, V_\ds(x) }{\ell(y) x^{1-\al \rho}} \,h_{\al \hat\rho}(y/x) } 
\quad &\mbox{as  $y\to\infty$ uniformly for  $1\leq x\leq y$}, \\[4mm]
{\displaystyle \frac{\al\hat\rho\, U_\as   (y)}{\hat\ell(x)y^{1-\al\hat\rho}} \, h_{\al\rho}(x/y )} \quad &\mbox{as  $x\to\infty$ uniformly for  $1\leq y \leq x$},
\end{array} \right.
\eeqn
\v2\nin
where 
$$h_\lam(\xi)  = \lam \int_0^1 t^{\lam-1} (\xi-1+t)^{\al-\lam-1}dt\quad (0<\lam \leq1, \xi\geq 1).$$

{\rm  (ii)}  Let $\al=1$ and $0<\rho<1$.  If $\rho\neq 1/2$,  then for each $0<\de<1$ the equivalence (\ref{asymp_g})  holds uniformly  both for   $1\leq x <\de y$ and for  $1\leq  y <\de x$,  and as $x\to\infty$
\beqn\label{eq(ii)} g_{\vom}(x,x) \sim \rho\hat\rho \int_0^x \frac{dt}{\ell(t)\hat\ell(t) t}; \eeqn
 and if  $F$ is recurrent, then
$$
g_\vom(x,x) \sim a(x) \sim \frac{2\sin^2 \rho\pi}{\pi^2} \int_0^x \frac{ dt}{L(t)t}.
$$
In either case,  for each $\ep >0$
\beqn\label{g/g1}
g_\vom(x,y) = o(g_\vom(y,y)) \quad  \mbox{as\; $y\to\infty$
 uniformly for\, $x: |x-y|> \ep y$}.
\eeqn

{\rm  (iii)}   If $F$ is transient, then $g_\vom(x,x) \to 1/P[\sigma_0=\infty]$ and 
for $x\geq 0$, $g_\vom(x,y) = o(g_\vom(y,y))$ as $|y-x|\wedge y \to \infty$. 
\end{lemma}
Note that $h_{\lambda}\equiv 1$ for $\lambda=\al-1$, $h_\lambda(1)=\lambda/(\alpha-1)$, and 
$$h_\lambda(\xi) \sim \xi^{\alpha-\lambda-1}\quad \mbox{as}\quad \xi\to\infty.$$
\pf\,  Let $1\leq x\leq y$. Then 
$$ g_{\vom}(x,y) =\sum_{k=0}^{x} v_\ds(k)u_\as   (y-x+k).$$
If $x/y\to 0$,  then 
$$ g_{\vom}(x,y) \sim \al\rho\, V_\ds(x)y^{\al\rho-1}/\ell(y),$$
which coincides with the asserted formula since $h_{\al\hat\rho}(\xi) \sim  \xi^{\al\rho-1}$ as  $\xi\to\infty$.
For $y\asymp x$ by Lemma \ref{lem*1}  the above sum divided by $\al^2\rho\hat\rho\,$  is asymptotically equivalent to
$$\sum_{k=0}^{x} \frac{k^{\al\hat \rho-1} (y-x+k)^{\al\rho-1}}{\hat\ell(k)\ell(y-x+k)}\sim \frac{x^{\al-1}}{\hat\ell(x)\ell(y)} \int_0^1 t^{\al\hat \rho-1} \bigg(\frac{y}{x}-1+t\bigg)^{\al\rho-1}dt
 \sim \frac{V_\ds(x)h_{\al\hat\rho}(y/x)}{\al\hat\rho\ell(y) x^{1-\al\rho}},$$
verifying the first formula of (\ref{asymp_g}). The second one is dealt with in the same way. 
 (i) has been proved. (iii) is easy to see (cf. Appendix (B)). 

Let $\al=1$ and  $0<\rho <1$. Then for $\rho\neq 1/2$,  by Lemma \ref{lem*1} $v_\ds(k) u_\as   (k) \sim \rho\hat\rho /[k\ell(k)\hat\ell(k)]$  and  (\ref{eq(ii)})  follows immediately. 
The first assertion of the case $\rho\neq \frac12$ is verified in the same way as for (i).
If  $F$ is recurrent, then  $a(x) \sim a(-x) \sim \pi^{-2} 2(\sin \rho\pi)^2\int_0^x[L(t)t]^{-1}dt$ according to \cite[Proposition 61(iv)]{Upot}.  Hence the remaining results  follow from Proposition \ref{prop3.1} (see Remark \ref{rem3.1}(c) and Lemma \ref{lem4.3}). \qed
\v2

Let $1<\al \leq 2$. Since $h_\lam(\xi) \sim \xi^{\al-\lam-1}$ as $\xi\to\infty$,  Lemma \ref{lem*2}(i) entails that
\beqn\label{asymp}
g_{\vom}(x,y) \asymp  \left\{ \begin{array} {ll} 
{\displaystyle V_\ds(x) U_\as   (y)/y, } 
\quad &\mbox{ for  $0 \leq x \leq y$},\\[2mm]
 V_\ds(x) U_\as   (y)/x \quad &\mbox{ for  $0 \leq y \leq  x$},
\end{array} \right.
\eeqn
where the constants involved in $\asymp$  depend only on $\al\rho$ and $\asymp$ can be replaced by $\sim$  in case $y/x\to \infty$ or $0$. By  Lemma 5.1(i) of \cite{Uexit} that gives the asymptotics  of $V_\ds(y)U_\as(y)H(y)$, it also follows that as $y\to \infty$
\beqn\label{asymp2}
g_{\vom}(y,y) \sim  \al\rho  h_{\al\hat\rho}(1)\frac{V_\ds(y)U_\as(y)}{y}  \sim  
\left\{\begin{array}{ll}
\lambda_{\al,\rho}/[yH(y)], \quad\; &1<\al<2,\\[2mm]
y/\int^y_0tH(t)dt,  \quad &\al=2,
\end{array}\right.
\eeqn
where  $\lambda_{\al,\rho} = \kappa \al^2\rho\hat\rho/(\al-1)>0$ ($\kappa$ is explicitly given as a function of $\rho$ and $\alpha$ only).


\begin{lemma} \label{lem7.3}  If $1<\al\leq 2$,  then as $R\to\infty$
\beqn\label{eqL3}
P_x\big[\sigma_R <T\big] \sim \left\{ \begin{array}{ll}
{\displaystyle  \bigg[ \frac{ (R/x)^{1-\al\rho}h_{\al\hat\rho}(R/x)}{h_{\al\hat\rho}(1)}\bigg] \frac{V_\ds(x)}{V_\ds(R)} }\quad \mbox{uniformly for} \; 1 \leq x\leq R,\\[6mm]
{\displaystyle  \bigg[ \frac{ (R/x)^{\al\hat\rho}h_{\al\rho}(x/R)}{h_{\al\rho}(1)}\bigg] \frac{V_\ds(x)}{V_\ds(R)} } \quad \mbox{uniformly for} \;  x\geq R;
  \end{array}\right.
 \eeqn
in particular
\beqn\label{eqL31}
P_x\big[\sigma_R <T\big] \left\{ 
\begin{array}{ll} \to 1 \quad &\mbox{as}\;\; x/R \to 1,\\
\sim  [(\al-1)/\al\rho] (R/x)^{1-\al\hat\rho} \hat\ell(R)/\hat\ell(x)\,  \quad &\mbox{as}\;\; x/R \to \infty,\\
\sim [(\al-1)/\al\hat\rho]\, V_\ds(x)/V_\ds(R) \quad &\mbox{as}\;\; x/R \to 0.
\end{array}\right.
\eeqn
\end{lemma}
\pf\,   Because of
the identity $P_x\big[\sigma_R <T\big] =g_{\vom}(x,R)/g_{\vom}(R,R)$
  the first formula  of (\ref{eqL3}) follows from Lemma \ref{lem*2}. 
  The derivation of the second one is similar.    By $\lim_{\xi\to\infty}\xi^{\al\rho-1}h_{\al\hat\rho}(\xi) = 1$  (\ref{eqL31})  follows from  (\ref{eqL3}) together with   $h_{\lam}(1) =\lam/(\al-1)$.
 \qed
 \v2
 
 The function $\xi^{1-\al\rho}h_{\al\hat\rho}(\xi)$ decreasingly approaches unity as $\xi\to\infty$ if $\al\rho<1$ and $h_{\al\hat\rho}(\xi) \equiv 1$ if $\al\rho=1$. Combining  (\ref{upp_bd}) with Lemma \ref{lem7.3}  yields that if $1<\al<2$, for $0\leq x\leq R$
\beqn\label{asymp3}
\frac{V_\ds(x)}{V_\ds(R)}\geq  P_x(\La_R)\geq P_x\big[\sigma_R <T\big]\geq \frac{\al-1}{\al\hat\rho}\cdot \frac{V_\ds(x)}{V_\ds(R)}\{1+o(1)\}.
\eeqn

{\bf Proof of  Proposition \ref{prop1.4}.} \, Let   $\alpha>1$. If $q>0$ we have $\al\hat\rho<1$   so that 
 $h_{\al\rho}(\xi) \sim \xi^{-\gamma}$ ($\xi\to\infty$)  with $\gamma =1-\alpha\hat\rho>0$,  and accordingly the second case of (\ref{eqL3})
implies that  for each $\ep>0$,
$$\liminf_{R\to\infty}\inf_{z\geq (1+\ep)R} P_z[T<\sigma_R] >0,$$ 
while if $p>0$ it follows that for a small $\ep>0$ and for all sufficiently large $R$
 \beqn\label{pr_(i)}
P_x\big[Z(R) > \ep R\,\big|\, \La_R\big] \geq 1/2\qquad  (x<\de R)
\eeqn 
 owing to Lemma 5.4 of \cite{Uexit}.
  These together shows that  the inequality $P_x\big[\sigma_R <T\big] \leq \theta V_\ds(x)/V_\ds(R)$  of (i)  holds under $pq> 0$.   
Now the  equivalence (\ref{Q_PLa})  in (i) 
follows from Lemma 47 of \cite{Upot}, that also shows that if $p=0$, $P_x[\sigma_R<T] \sim P_x(\La_R)$ uniformly
for $0\leq x\leq R$ (this may be verified directly because of  the first case of (\ref{eqL31})). 
If  $q=0$, then $\alpha\hat \rho= 1$ so that $h_{\al\rho} =h_{\al-1} \equiv 1$ and (\ref{eqL3}) implies that  $P_z[\sigma_R<T]\to 1$ uniformly for  $R<z<MR$ for any $M>1$ and we conclude that $P_x[\sigma_R<T]\sim P_x(\La_R)$ uniformly for $0\leq x\leq R$ since for such $x$,  $P_x[Z(R)<MR\,|\, \La_R]\to 1$ as $R\to\infty$ and $M\to\infty$.
 The asymptotic equivalence stated last in (i) is a reduced form of the first  formula  in  (\ref{eqL3}).




(ii) follows  from  (\ref{g/V/a}) of  Theorem \ref{thm1} 
with the help of (\ref{fun_res}) and the asymptotic forms of $a(x)$ and $a(-x)$ given in (\ref{P_a})).

The case $y\to\infty$ of (iii) follows from (\ref{T2-3}) of Theorem  \ref{thm2}. The other case is cheaper and immediate from (\ref{Dual}). 

  (iv)   follows from Lemma \ref{lem7.4} given in the next subsection if $\rho=1$.
 In the other case    $0<\rho<1$,  with the help of  Lemma 5.4 of \cite{Uexit} that says that if $0<(\al\vee1)\rho<1$ we have 
 \beqn\label{Ovsh}
 P_x\big[Z(R)\leq \ep R\,\big|\, \La_R\big] \to 0  \quad \mbox{as \;\; $R\to\infty$ and $\ep\downarrow 0$}
\eeqn
uniformly for $0\leq x <\de R$, and the required convergence follows.

 \v2\v2
 
{\bf 7.2.}  {\sc Transient walks.} 
\v2
   In  \cite{Uexit} we have brought in the condition
\v2
(C4)  \quad  (AS)  holds with  \, $\alpha<1=\rho.$
\v2\noindent
 If  either (C3) or (C4)  holds, then  $u_\as   (x)\sim 1/\ell^*(x)$ so that 
\beqn\label{4.15}
P_x[\sigma_R<T]=\frac{g_{\vom}(x,R)} {g_{\vom}(R,R)} \sim \frac{V_\ds(x)U_\as   (R)}{Rg_\vom(R,R)} \quad \mbox{uniformly for} \; 0\leq x<\de R,
\eeqn
whether $F$ is recurrent or transient.  
  In the next lemma  we show what is asserted in (iv) of Proposition \ref{prop1.4} in case $\rho=1$, when either (C3) or (C4)  holds.   
  Under (C4)  $ \ell_\sharp$ should be defined by 
 \beqn\label{el/sh}
  \ell_\sharp(t) = \alpha \int_t^\infty \frac{s^{\al-1}F(-s)}{\ell (s)}ds \quad (t >0).
  \eeqn
 \begin{lemma}\label{lem7.4}\, Let $F$ be transient.  Then 
 \beqn\label{4.16}
 g_\vom(x,x) \to 1/q_{\infty} \,(<\infty),
 \eeqn
and  if either (C3) or (C4) holds, for each  $\de<1$, as $R\to\infty$
  \beqn\label{4.17}
  P_x[\sigma_R <T\,|\, \La_R]\sim q_\infty/\ell^*(R)\ell_\sharp(R) \,  \longrightarrow\, 0\quad \mbox{ uniformly for}\; \;  0\leq x\leq \de R.
  \eeqn
\end{lemma} 
\pf 
(\ref{4.16}) is a standard result for a  general transient r.w.\,(cf. Appendix (B)).
The equivalence  in (\ref{4.17})  follows from (\ref{4.16}) in view of (\ref{4.15}). If (C4) holds, this entails (\ref{4.17}), for by Lemma \ref{lem6.4} of Appendix (A)  $\ell^*(x)\ell_\sharp(x) \sim A(x) \to\infty$.

 Suppose (C3) to hold. By the transience of $F$---entailing $E|X|=\infty$---the probability $P_x[\sigma_0<\infty]$ tends to zero, hence  $P_{x+w}[\sigma_R <\infty]\to 0$ as $w\to\infty$, and it accordingly suffices to show that for any constant $M>1$,
\beqn\label{4.18}
P_x[Z(R) <M\,|\,\La_R] \to 0 \quad (R\to\infty)
\eeqn
uniformly for $0\leq x<\de R$.  Put $B= (-\infty,-1]\cup [R+1,\infty)$. 
Then the conditional probability above is expressed as
$$\frac1{P_x(\La_R)}\sum_{w=0}^R g_B(x,R-w) P[0 < X-w<M].
$$
Let $N= \lfloor (1-\de)R/2\rfloor$ so that $R-x \geq 2N$ for $0\leq x\leq \de R$. We claim that
\beqn\label{4.19}
\sum_{0\leq w\leq N} g_{[R,\infty)}(x,R-w) P[0 < X-w<M] = o(P_x(\La_R))\quad (x\leq \de R).
\eeqn
Since $g_B(x,R-w) \leq g_{[1,\infty)}(x-R,-w) = g_{\vom}(w, R-x) \sim V_\ds(w)/\ell^*(R)$ for  $1\leq w\leq N$ and,  since $E[V_\ds(X); X\geq 0] = V_\ds(0)$, we obtain
\beqn\label{4.20}
\sum_{0\leq w\leq N} g_{B}(x,R-w) P[X-w=y]\leq C/\ell^*(R) \quad (0\leq y\leq M).
\eeqn
Summing over $y$ yields that the sum on the LHS is at most a constant multiple of $M/\ell^*(R)$ which is $o(P_x(\La_R))$ as $x\to\infty$, entailing (\ref{4.19}), for  
$V_\ds(R)/\ell^*(R) \sim 1/[\ell^*(R)\ell_\sharp(R)]$ is bounded owing to the equivalence in (\ref{4.17}) that we have already seen to be true.  When $x$ remains  in a bounded interval,
(\ref{4.19})  also follows from (\ref{4.20}). Indeed, for each $x$ fixed, taking any $0<\ep<1/2$ and choosing  a constant $r=r(\ep,x)$  so that $P_x(\La_r)<\ep$ yield that for $0\leq w\leq N$,
$$g_B(x, R-w) = \sum_{1\leq z\leq R/2} P_x[ Z(r) =z, \La_r]g_B(r+z, R-w) + o(P_x(\La_R))$$
for $P_x[Z(r) > R/2\,|\,\La_r] \leq Cr H_+(R/2) = o(P_x(\La_R))$, and substituting this into (\ref{4.20}) and taking summation over $w$ first, as in the same way as above we have
$$\sum_{0\leq w\leq N} g_{B}(x,R-w) P[X-w=y]\leq C\ep/\ell^*(R)+o(P_x(\La_R)) \quad (0\leq y\leq M).
$$
Hence we obtain (\ref{4.19}), $\ep$ being arbitrary. By $g_B(x,x)\leq C_1$ it is easy to see
$$\sum_{N< w \leq R} g_{B}(x,R-w) P[0< X-w<M]\leq C_1MH_+(N) =o(P_x(\La_R)),
$$
which together with (\ref{4.19}) concludes (\ref{4.18}). Proof of Lemma \ref{lem7.4} is finished. 
 \qed 
 
\section{Appendix}

\quad {\bf (A)} \;  
  Let $\ell_\sharp$ be  the function defined by (\ref{3.0}) in the preceding section. The following result is taken from \cite{Uexit}. The condition (C4) is given at  the beginning of Section 7.2 and 
 the definition of $\ell_\sharp$ must be given by (\ref{el/sh}) rather than (\ref{3.01}). 
\begin{lemma}\label{lem6.4}  \, Suppose that either (C3) or (C4) holds. Then
 \beqn\label{id_A}
\ell^*(t) \ell_\sharp(t)= -\int_0^t F(-s)ds +\int_0^t P[Z>s]\ell_\sharp(s)ds
=A(t) + o\bigg(\int_0^t [1-F(s)]ds \bigg)
 \eeqn
and in case $EX=0$,  both $\eta_-$ and $\eta$ are s.v. and
\beqn\label{id_A2}
\ell^*(t) \ell_\sharp(t) = \int_t^\infty\big[F(-s) - P[Z>s]\ell_\sharp(s)\big]ds=A(t) +o(\eta_+(t)).\eeqn
 \end{lemma}
 
 \v2
 Put $A_\pm(x) = \int_0^x P[\pm X>t]dt$ and suppose (C3) to hold. If the positive and negative tails of $F$ are not balanced in the sense that
 \beqn\label{el/A}
 \limsup \frac{\eta_+(x)}{\eta_-(x)} <1 \quad \mbox{if}\;\; EX=0\quad \mbox{and}\quad \limsup \frac{A_-(x)}{A_+(x)} <1\quad \mbox{if}\;\; E|X|=\infty,
 \eeqn
 then  (\ref{id_A}) and (\ref{id_A2})
 together show (\ref{A/el/a}), i.e., $ \ell^*(x)\ell_\sharp(x) \sim A(x)$,
  or equivalently in view of   {\bf  L(3.1)},
 \beqn\label{eqR}
 V_\ds(x)U_\as   (x) \sim x/A(x).
  \eeqn 
  Since $\eta$ is s.v. (as noted in Lemma \ref{lem6.4}), hence $xH(x)/\eta(x) \to 0$ and, under the first condition of (\ref{el/A}),  $A(x) =\eta_-(x)-\eta_+(x) \asymp \eta_-(x) \asymp  \eta(x)$, it follows that $F$ is p.r.s.; in particular  if $F$ is recurrent,
$a(x) \asymp \int_{0}^x \big[ F(-t)/\eta_-^2(t)\big]dt \sim 1/\eta_-(x)$, hence $a(x)\asymp 1/A(x)$, or what amounts to the same,  $\limsup a(-x)/a(x)<1$. Thus if $EX=0$, 
\beqn\label{a/eta}
 \limsup \frac{\eta_+(x)}{\eta_-(x)} <1 \,\Longrightarrow \, 
  \limsup \frac{a(-x)}{a(x)} <1 \,\Longrightarrow \, \mbox{(\ref{eqR})}.
\eeqn 
where the second implication is observed in Remark \ref{lem3.1}(b). We do not know whether the converse of the first implication in (\ref{a/eta}) is true or not. 
 From (\ref{id_A2}), 
 being written as $\ell^*(x)\ell_\sharp(x) =\eta_-(x)-\eta_+(x)+o(\eta_+(x))$, we  also infer  that  under $EX=0$,  
$$\limsup \frac{\eta_-(x)}{A(x)}<\infty \, \Longleftrightarrow \, \limsup\frac{\eta_+(x)}{\eta_-(x)}<1  \, \Longleftrightarrow \, \limsup \frac{\eta_-(x)}{\ell^*(x)\ell_\sharp(x)} <\infty,$$
and combining this with
(\ref{aaa}) we see that  $\limsup a(x)\eta(x) =\infty$ if $\limsup a(-x)/a(x)=1$.  
  
Suppose that $F$ is transient and p.r.s..  Then by the dual of 
Lemma \ref{lem5.1} (see (a) right after it) $g_\vom(x,2x) = G(x)-G(-x) +o(G(x))$, and the same reasoning as for the recurrent $F$ shows that   $\limsup G(-x)/G(x)<1$ implies (\ref{eqR}).   In a similar way to the recurrent case, using $G(x)\sim \int_x^\infty \big[H_+(t)/A^2(t)\big]dt$, 
we see that 
\beqn\label{A/G}
 \limsup \frac{A_-(x)}{A_+(x)} <1 \,\Longrightarrow \, 
  \limsup \frac{G(-x)}{G(x)} <1 \,\Longrightarrow \, \mbox{(\ref{eqR})}.
\eeqn

\v2
{\bf (B)}  \; Let $F$ be transient so that we have the Green kernel 
$G(x):= \sum_{n=0}^\infty P[S_n=x]<\infty$.
For  $y\geq 0, x \in \mathbb{Z}$,
 \[
 G(y-x) -g_{\vom}(x,y) = \sum_{w=1}^\infty P_x[S_T=-w]G(y+w).
 \]
 According to the Feller-Orey renewal theorem \cite[Section XI.9]{F}, $\lim_{|x| \to \infty} G(x)=0$ (under $E|X|=\infty$),  showing that the RHS above tends to zero as $y\to\infty$ (uniformly in $x\in \mathbb{Z}$), in particular $\lim g_{\vom}(x,x) =G(0) =1/P[\sigma_0=\infty]$. It also follows that  $P_x[\sigma_0<\infty]= G(-x)/G(0) \to0$.
 
 \v2
{\bf (C)} \;  Let $T_0=0$ and $T= (T_n)_{n=0}^\infty$ be a r.w. on  $\{0,1,2,\ldots\}$  with  i.i.d.   increments. Put
$$u(x) = \sum_{n=0}^\infty P[T_n=x]\;\;  (x = 1,2,\ldots). $$
Suppose that $T_1$ is aperiodic so that  $u(x)$ is positive for all sufficiently large $x$. 
We give some results as to  asymptotics of $u(x)$ when
 the tail  $$\ell(t):=P[T_1>t]$$
  is s.v., or what is the same thing the renewal function
  $$U(x)= 1+u(1)+\cdots + u(x)$$
  is s.v.,
    Nagaev \cite{Nv} shows that if  $xP[T_1=x]$ is s.v., then $u(x)\sim P[T_1=x]/[\ell(x)]^2$.   
  For the proof of Theorem \ref{thm3} we need the upper estimate $u(x)= o(U(x)/x)$.  The following lemma,  slightly extending \cite{Nv},  gives a better bound under a restriction on 
  $$q(x):= P[T_1=x].$$
 \begin{lemma}\label{lem6.2.2} Suppose  $P[T_1>t]$ is s.v.  If for a  constant $C$  and $x_0$
 \beqn\label{eL2.2}  C:=\lim_{\de\uparrow 1} \limsup_{x\to\infty} \frac1{q(x)} \sup_{\de x< y\leq x}  q(y) <\infty,
 \eeqn 
 then
 $$u(x) \leq \frac{q(x)}{[\ell(x)]^2}\{C+o(1)\}.$$
If 
\beqn\label{eL2.2*}   c:=\lim_{\de\uparrow 1} \liminf_{x\to\infty} \frac1{q(x)} \inf_{\de x< y\leq x}  q(y) <\infty,
  \eeqn
  then 
  \[
  u(x) \geq \frac{q(x)}{[\ell(x)]^2}\{c+o(1)\}.
  \]
In particular if both (\ref{eL2.2}) and (\ref{eL2.2*}) holds with  $C=c=1$, then  $u(x) \sim q(x)/[\ell(x)]^2.$
 \end{lemma}
 \v2
 \pf The proof  is made by elaborating on that of  \cite{Nv}.
 Put
$$u^{(2)}(x)= \sum_{y=0}^x u(x-y)u(y).$$
Then it holds  \cite[Lemma 2.6]{Nv}  that 
\beqn\label{eL2.21}
xu(x) =\sum_{y=0}^{x-1} (x-y) q(x-y) u^{(2)}(y).
\eeqn
This is derived by means of the generating functions. Indeed 
if $f(s)=\sum q(x)s^z$  and  $h(s) =  \sum u(x)s^x$ ($|s|\leq 1$)  (the generating functions of  $q(\cdot)$ and $u(\cdot)$), then $h(s)=1/[1-f(s)]$.    The identity (\ref{eL2.21}) follows by comparing  the identities  $\sum (x+1)q(x+1)s^x =f'(s)$, 
$$\sum u^{(2)}(x)s^x = \frac1{[1-f(s)]^2}\quad\mbox{and}\quad
 \sum (x+1)u^{(2)}(x+1)s^x = h'(s) \quad (|s|<1).$$  


The slow variation of $\ell$ entails
\beqn\label{eL2.25} U(x)\sim 1/\ell(x),\quad\quad u^{(2)}(x) \leq  \frac{2+o(1)}{\ell(x)}\max_{\frac12 x\leq y\leq x} u(y),
\eeqn
and
$$\sum_{y=0}^x  u^{(2)}(y) =\sum_{y=0}^x u(y)\sum_{z=y}^x u(z-y) = [U({\textstyle \frac12 }x)]^2\{1+o(1)\}+ \sum_{y=x/2}^{x}U(x-y)u(y).$$
The last sum is less than  $U(x) \sum_{y=x/2}^x u(y) = o([U(x)]^2)$ it follows that
\beqn\label{eL2.23}
\sum_{y=0}^x  u^{(2)}(y)  \sim \frac1{[\ell(x)]^2}.
\eeqn 

Let $0<\de<1$ and $\ep=1-\de$, put 
$$C_\de =\de^{-1}\limsup_{x\to\infty}\sup_{\de x\leq y\leq x} q(y)/q(x)$$
and split the  range of the sum on the RHS of (\ref{eL2.21})  according as  $y\leq \ep x$ or $y>\ep x$. Then using (\ref{eL2.25} ) and (\ref{eL2.2}) one sees that for all sufficiently large  $x$,
\beqn\label{eL2.22}
xu(x) \leq C_\de \frac{xq(x)}{[\ell(x)]^2} + \frac{3\sum_{y=0}^{\de x} yq(y)}{\ell(x)}\max_{\ep x\leq y\leq x} u(y).
\eeqn
Since $\sum_{y=0}^x yq(y)= - x\ell(x) +\sum_{y=0}^{x-1} \ell(y) =  o(x\ell(x))$,  on writing 
$$ N(x) = \frac{q(x)}{[\ell(x)]^2},  \quad  M_\de(x) = \max_{\ep x\leq y\leq x} u(y)$$
 this yields that for $x$ large enough
$$u(x) \leq C_\de N(x) + o(M_\de(x)).$$

Condition (\ref{eL2.2}) entails that $q(x)\geq c_1q(y)$ for $x/2\leq y< x$ with  $c_1>0$, and accordingly   one can find  positive constants $R_0$ and $\alpha$ such that $q(x) > x^{-\alpha}$ for $x\geq R_0$.   Take $\eta=\eta_\de >0$ such  that 
$$\eta C_\de< \ep \quad\mbox{and}\quad   \log \eta^{-1} > 2\alpha \log \ep^{-1}.$$
 Choose $R\geq R_0$ so that  for $x>R$
$$u(x) \leq  C_\de N(x) + \eta M_\de(x) \quad\mbox{and}\quad \sup_{\de x\leq y<x} N(y)<  C_\de N(x).$$
 Take $x_1\in [\ep x, x]$ such that  $u(x_1)=M_\de(x)$ so that
$$u(x) \leq C_\de N(x) + \eta  u(x_1).$$ 
If $x\geq R/\ep$, then $u(x_1) \leq  C_\de N(x_1) + \eta M_\de(x_1)$, hence 
$$u(x) \leq C_\de \big(N(x) + \eta  N(x_1) \big)+ \eta^2 M_\de (x_1).$$
If $x\geq \ep^{-k-1}R$ one can repeat this procedure $k$ times  to obtain 
$$u(x) \leq C_\de \big[N(x) + \eta N(x_1) +  \cdots + \eta^k N(x_k) \big] + \eta^{k+1} M_\de(x_k).$$
Let  $n(x)$ be the largest integer  $n$ such that $x\geq \ep^{-n-1}R$.  Since    $N(x_j) \leq  C_\de^jN(x)$,  recalling $\eta  C_\de \leq \ep$  one infers that
$$u(x) \leq \de^{-1}C_\de N(x) + \eta^{n(x)} M_\de(x_{n(x)}).$$
One can easily see that $\eta^{n(x)}/N(x)\to 0$. Hence 
$\limsup u(x)/N(x) \leq \de^{-1}C_\de$. This concludes the asserted upper bound, for $\de^{-1}C_\de$ can be made  arbitrarily close to $C$.

The lower bound is easily deduced from (\ref{eL2.21}). Indeed, the restriction to $y\leq (1-\de)x$ of the sum on its RHS is larger than
$$(1-\de)x \inf_{\de x\leq y\leq x} q(y) \sum_{y=0}^{(1-\de)x}u^{(2)}(y),$$
and (\ref{eL2.2*}) together with (\ref{eL2.23})  yields the asserted lower bound of $u(x)$.  \qed

\end{document}